\documentclass[a4paper,10pt,twocolumn]{article}
\usepackage[utf8]{inputenc}

\usepackage[switch]{lineno}

\usepackage{hyperref}
\usepackage{bm}
\usepackage{authblk}
\usepackage{siunitx}
\usepackage{graphicx}
\usepackage{amsmath, amssymb}
\usepackage{subcaption}  

\usepackage[a4paper,left=1.5cm,right=1.5cm,top=2cm,bottom=2cm]{geometry}

\usepackage{nomencl}
\usepackage{ifthen}
\renewcommand{\nomgroup}[1]{%
\ifthenelse{\equal{#1}{C}}{\item[\textbf{Acronyms}]}{%
\ifthenelse{\equal{#1}{G}}{\item[\textbf{Greek Letters}]}{%
\ifthenelse{\equal{#1}{L}}{\item[\textbf{Latin Symbols}]}{%
\ifthenelse{\equal{#1}{S}}{\item[\textbf{Subscripts}]}{}}}}
}
\makenomenclature

\usepackage[numbers,sort&compress]{natbib}

\def\tsc#1{\csdef{#1}{\textsc{\lowercase{#1}}\xspace}}
\tsc{WGM}
\tsc{QE}
\tsc{EP}
\tsc{PMS}
\tsc{BEC}
\tsc{DE}

\begin{document}
\let\WriteBookmarks\relax
\def\floatpagepagefraction{1}
\def\textpagefraction{.001}

\title{Deep Reinforcement Learning for the Heat Transfer Control of Pulsating Impinging Jets}
\author[1]{Sajad~Salavatidezfouli\footnote{ssalavat@sissa.it}}
\author[2]{Giovanni~Stabile\footnote{gstabile@sissa.it}}
\author[1]{Gianluigi~Rozza\footnote{grozza@sissa.it}}

\affil[1]{Mathematics Area, MathLab, International School for Advanced Studies (SISSA), Trieste, Italy}

\affil[2]{Department of Pure and Applied Sciences, Informatics and Mathematics Section, University of Urbino Carlo Bo, Urbino, Italy}

\date{} 

\twocolumn[
  \begin{@twocolumnfalse}
    \maketitle
	   \begin{abstract}
        This research study explores the applicability of Deep Reinforcement Learning (DRL) for thermal control based on Computational Fluid Dynamics. To accomplish that, the forced convection on a hot plate prone to a pulsating cooling jet with variable velocity has been investigated. We begin with evaluating the efficiency and viability of a \textit{vanilla} Deep Q-Network (DQN) method for thermal control. Subsequently, a comprehensive comparison between different variants of DRL is conducted. Soft Double and Duel DQN achieved better thermal control performance among all the variants due to their efficient learning and action prioritization capabilities. Results demonstrate that the soft Double DQN outperforms the hard Double DQN. Moreover, soft Double and Duel can maintain the temperature in the desired threshold for more than 98\% of the control cycle. These findings demonstrate the promising potential of DRL in effectively addressing thermal control systems.

    	\vspace{0.5cm}
		\textbf{Highlights:}
		\begin{itemize}
		  \item DRL technique demonstrates a successful thermal control performance.		
		  \item Hard Double DQN is not useful for thermal control tasks.
		  \item Soft Double and Duel DQN demonstrate more temperature uniformity along the surface.
    
		\end{itemize}

		\end{abstract}
  \end{@twocolumnfalse}

]
\maketitle
\printnomenclature

\section{Introduction}
Thermal control is an essential practice with a profound impact on diverse applications, spanning HVAC systems, electronics cooling, medical devices, food and beverage production, and data centers. Achieving optimal operation in these fields often depends on the maintenance of the temperature within a stable and narrow threshold. This can be addressed by the manipulation of heat transfer mechanisms, i.e. conduction, convection, and radiation \citep{bergman2011introduction, childs1999heat}. Among these mechanisms, convection holds a prominent position, where the transfer of heat is facilitated through the movement of surrounding fluid \citep{bejan2013convection}.

In recent years, convection control gained a lot of concentration \citep{giwa2021review, xiong2021numerical, shahabadi2021controlling}. Notably, forced convection control has emerged as a focal point, owing to its inherent advantages. Forced convection accelerates the rate of heat transfer, making it a more effective method that can better meet temperature requirements. Pioneering studies have shed light on the potential of this approach. For instance, Davalath and Bayazitoglu examined the forced convection between parallel plates consisting of finite block heat sources \citep{davalath1987forced}. They could regulate the flow and temperature fields through the spacing of blocks, Reynolds and Prandtl numbers. Al-Sarkhi and E. Abu-Nada focused on thermal control optimization by modifying the height and number of fins in a tube \citep{al2005characteristics}, pinpointing an optimal combination of these control parameters. In a computational study undertaken by Yilmaz and Oztop, turbulent forced convection in a double forward-facing step was investigated \citep{oztop2006turbulence}. Their findings underscored the significance of step size as a passive control element in heat transfer scenarios. Extending this pursuit of enhanced heat transfer, Kim et al. employed a control strategy by means of a synthetic jet in a channel \citep{kim2022enhancing}. Their results proved the superiority of higher frequencies and direct impingement in achieving enhanced cooling performance. All of these studies exclusively employed passive thermal control techniques.

Numerous studies have explored active thermal control to enhance the forced convection mechanism \citep{tan2021general, xu2021near, zhang2019switched}. Active control systems, which rely on external power for operation, have shown particular effectiveness in achieving a precise temperature range. Notably, these investigations have highlighted the potential of active methods to deliver highly targeted thermal control systems \citep{kuprat2021research}.

However, it is essential to acknowledge that many of these studies require intensive numerical resources to derive their control strategies \cite{miao2021spacecraft}. This reliance on computational resources has raised concerns about the practicality and efficiency of these approaches in modern industrial settings. Furthermore, while active control systems are promising, there remains a significant challenge in directly applying optimal control algorithms due to the inherent complexity of governing partial differential equations (PDEs). Despite ongoing research, there is no universally applicable and efficient method for convection control that has emerged as of yet \citep{kuprat2021research}.

Over the past few years, the expansion of high-performance computing resources has paved the way for the integration of data-driven and non-intrusive machine learning (ML) methods into the field of fluid dynamics \citep{taira2017modal, brenner2019perspective, brunton2020machine, taira2020modal}. ML's adaptability offers an innovative modelling framework capable of effectively addressing numerous challenges within fluid mechanics. These challenges encompass reduced-order modeling \citep{bai2017data, san2019artificial, peng2020time, peng2021geometry}, shape optimization \citep{li2020efficient}, and turbulence closure modeling \citep{pawar2020priori}. Moreover, ML-based approaches have demonstrated significant promise in heat transfer problems, including contact heat transfer \citep{vu2021machine}, critical heat flux \citep{swartz2021machine}, and thermal resistance \citep{wu2019predicting} prediction and optimization.

Active thermal control systems encounter a multitude of challenges, ranging from achieving precision and energy efficiency to addressing environmental concerns and controlling heat dissipation. However, one important challenge that ties these together is the time delay. Central to any thermal control system is the time required for heat transfer to occur, allowing the controller to make informed decisions. Fortunately, deep reinforcement learning (DRL), a novel branch of machine learning, has showcased its capacity to tackle both nonlinear challenges and time delay issues in flow control \citep{zou2020towards}.

In DRL problems, the choice of environment plays a crucial role in determining the success of control applications. An environment in DRL represents the simulated or real-world system with which an agent interacts to learn optimal control strategies. Various environments can employed in DRL, ranging from simple, 0-Dimensional models to more complex, 3-Dimensioanl simulations or even experimental data \citep{gao2020deepcomfort, wang2020reinforcement, xiong2022application}.

When it comes to controlling forced convection systems, the choice of environment becomes particularly critical. Many previous studies have explored DRL in conjunction with simplified, analytical models, which provide rough estimations of fluid behaviour. However, the precision of these simplified models may not always meet the demands of highly accurate thermal control. In contrast, Computational Fluid Dynamics (CFD) stands out as a superior choice for the environment. CFD is known for its ability to the complete form of the fluid flow equations, providing a highly accurate representation of fluid behaviour. 

While DRL-CFD studies remain relatively scarce in comparison to other DRL environments, they have demonstrated significant promise. These studies have progressed from computationally efficient investigations into laminar flow \citep{rabault2019artificial, tang2020robust, xu2020active, li2022reinforcement} to experimental research involving high Reynolds flow \citep{ren2021applying, shimomura2020experimental, fan2020reinforcement}. Furthermore, DRL-CFD has been successfully applied in areas such as flow separation suppression \citep{wang2022deep} and the enhancement of vortex-induced vibration \citep{mei2021active}. This emerging field holds great potential for effectively controlling forced convection systems, thanks to the precision and accuracy of CFD as the chosen environment.

While DRL-CFD studies remain relatively rare, their potential applications in forced convection control are highly promising. Notably, this study explores active thermal control by means of impingement cooling, as there is no similar study in the literature, to the best of our knowledge. Moreover, it is worth emphasizing that there is a notable absence of comparable research that systematically compares different DRL methods in the context of thermal control.

The effectiveness of the DRL-CFD is demonstrated by evaluating its cooling effect on a system comprising a heater surface under constant heat flux and a cooling jet with controlled velocity. In Section 2, we provide a detailed discussion of the methodology, covering topics such as the reinforcement learning framework, CFD solver, deep reinforcement learning, Deep Q-Networks (DQNs), and the algorithms employed in this study. Section 3 provides the model description. Finally, Section 4 focuses on the results of the implemented control system on the hot plate followed by the discussion on the findings.

\section{Reinforcement Learning}
Reinforcement learning (RL) is a fascinating branch of machine learning that utilizes a closed-loop feedback control framework. It represents a novel approach to automatically discovering the optimal control strategy. RL comprises various components and follows a well-defined execution process, depicted in Figure 1.

RL begins with the acquisition of an observation, referred to as the state ($s_t$), which is sampled from the environment at specific times. Based on this observation, the agent selects an action ($a_t$) that maximizes a specific value called the reward ($r_t$), which is calculated by a certain \textit{reward} function. The chosen action acts as a control signal that is executed within the environment, leading to the sampling of the next observation ($s_{t+1}$) at the subsequent time. This iterative interaction continues until a predefined terminal condition is met, such as environmental convergence or a designated duration, constituting an episode. Through training over numerous episodes, the agent gains the ability to excel in its performance. In other words, it becomes capable of generating an optimal trajectory of action-state pairs ($s_0$, $a_0$, $s_1$, $a_1$, ...).

There exist two primary branches of RL known as model-based and model-free \citep{doll2012ubiquity}. In the model-based method, the agent constructs a model during the training phase to capture the relationship between states and actions in the environment \citep{kaiser2019model}. This enables the agent to predict the outcome of executing any action in the environment. Numerous model-based algorithms have been developed particularly for solving kinetics and motion planning problems \citep{polydoros2017survey}. However, one noticeable limitation of model-based methods is their heavy reliance on a comprehensive understanding of the rules governing the environment. This requirement poses challenges, especially when dealing with nonlinear problems where constructing an accurate environment model becomes exceedingly difficult. In such cases, model-free methods can offer advantages over model-based approaches. Model-free methods do not necessitate an environment model of any kind, and an agent utilizing model-free methods explicitly learns through trial and error, gradually transitioning from tentative exploration to deliberate planning at higher levels of proficiency \citep{ccalicsir2019model}.

This research primarily concentrates on addressing a fluid dynamics problem namely the control of turbulent, incompressible fluid flow in conjunction with heat transfer. The problem at hand pertains to a 3-dimensional scenario, where the conservation of mass, momentum, and energy are governed by the interdependent and nonlinear Navier-Stokes and heat equations:
\begin{equation}
\label{eq1}
\nabla \cdot \boldsymbol{u}=0, \\
\end{equation}
\begin{equation}
\label{eq1_2}
\frac{\partial \boldsymbol{u}}{\partial t}+\boldsymbol{u} \cdot(\nabla \boldsymbol{u})=-\nabla p+\nu \nabla \cdot \nabla \boldsymbol{u}, \\
\end{equation}
\begin{equation}
\label{eq1_3}
\frac{\partial T}{\partial t}+\boldsymbol{u} \cdot \nabla T=\alpha \nabla \cdot \nabla T ,
\end{equation}
where, the velocity vector is represented by $\boldsymbol{u}$, pressure by $p$, temperature by $T$, kinematic coefficient of viscosity by $\nu$, and thermal diffusivity by $\alpha$. The numerical simulation is conducted with the assistance of ANSYS FLUENT \citep{ansys2019ansys}. It is worth noting that we have opted for a model-free approach, considering the inherent complexity of solving nonlinear partial differential equations.

\subsection{Deep Q Network}\label{section2_1}
The interaction between agent and environment can be formalized as a Markov Decision Process (MDP) in a typical reinforcement learning by a tuple ($S, A, T, r, \gamma$) where $S$ and $A$ are set of states and actions. $T\left(s, a, s^{\prime}\right)=P\left[S_{t+1}=s^{\prime} \mid S_t=s, A_t=a\right]$ is the transition function and reward function, $r$, is $r(s, a)=\mathbb{E}\left[R_{t+1} \mid S_t=s, A_t=a\right]$. Finally, $\gamma$ is the discount factor and is $\gamma \in[0,1]$. MDP is dedicated to the classical sequential decision-making process, where actions affect both the immediate reward and subsequent states. Thus, the discounted return is defined as $G_{t}=\sum_{k=0}^{\infty} \gamma_t^{(k)} R_{t+k+1}$ from the state as the discounted sum of ongoing rewards and $\gamma_t^{(k)}=\prod_{i=1}^k \gamma_{t+i}$. The main goal of the agent is to maximize the discounted return $G_t$.

Q-learning is a well-established and effective method \citep{jang2019q} that does not require prior knowledge of the system's dynamics, i.e. it is a model-free method. It allows the agent to learn an approximation of the expected discounted return, also known as the value, by formulating an action-value function as follows:

\begin{equation}
\label{eq2}
Q^\pi(s, a)=\mathbb{E}_\pi\left[G_t \mid S_t=s, A_t=a\right].
\end{equation}

Where expectation $\mathbb{E}_\pi$ represents the average value when the agent selects an action according to the policy $\pi$. As mentioned before, the main goal of reinforcement learning (RL) is to find a policy that maximizes the reward. There exists at least one optimal policy, $\pi^*$, which outperforms any other policy and achieves the optimal action-value function. Thus, we can express the optimal action-value function $Q^*$ as follows:

\begin{equation}
\label{eq3}
\begin{aligned}
& Q^*(s, a)=\max _\pi Q^\pi(s, a).
\end{aligned}
\end{equation}

To derive a new policy from the action-value function, a commonly used approach is $\epsilon$-greedy. Accordingly, the agent selects the action with the highest value (referred to as the \textit{greedy} action) with a probability of (1-$\epsilon$), while choosing a random action uniformly with a probability of $\epsilon$.

Q-learning simply aims to estimate the optimal action-value function, $Q^*$, through an iterative learning process. The key concept behind Q-learning is to update the action-value function, Q, by iteratively refining its estimates:

\begin{equation}
\label{eq4}
\begin{aligned}
Q_n(s, a)=Q_{n-1}(s, a)+\alpha & [r+\gamma \max _{a^{\prime}} Q_{n-1}\left(s^{\prime}, a^{\prime}\right) \\
& \left.-Q_{n-1}(s, a)\right].
\end{aligned}
\end{equation}

The foundation of Q-learning lies in an important identity known as the Bellman equation. Accordingly, if we have knowledge of the optimal value at the next timestep, $Q^{*}(s', a')$, for all possible actions, $a'$, at the next time step, then the best strategy is to select the action, $a'$, that maximizes the expected value of the sum of immediate reward, $r$, and the discounted future value, $\gamma Q^{*}(s', a')$:

\begin{equation}
\label{eq5}
\begin{aligned}
& Q^*(s, a)=\mathbb{E}_{s^{\prime}}\left[r+\gamma \max _{a^{\prime}} Q^*\left(s^{\prime}, a^{\prime}\right) \mid s, a\right].
\end{aligned}
\end{equation}

By using the Bellman equation in an iterative manner, we converge towards the optimal action-value function, $Qn \rightarrow Q^*$, as $n \rightarrow \infty$. Interestingly, the learned action-value function, $Q_n$, directly approximates $Q^*$ regardless of the policy being followed. This simplifies the analysis of the algorithm and facilitates early convergence proofs. Consequently, as long as all state-action pairs are continuously updated, the estimate will ultimately converge to the correct values.

However, dealing with problems that involve large states and/or action spaces, such as active flow control, becomes exceedingly challenging when attempting to learn Q-value estimation for all possible state-action pair. To address this issue, novel approaches have emerged, employing deep neural networks to represent different aspects of an agent including policies $\pi(s, a)$ and values $Q(s, a)$. Within these methods, deep neural networks are utilized as nonlinear function approximators, trained through gradient descent to minimize a suitable loss function and achieve accurate estimation. An exemplary example of this approach is the Deep Q-Network (DQN) proposed by Mnih et al. \citep{mnih2013playing, mnih2015human} in which the deep networks and reinforcement learning were successfully coupled. At each step, the agent, based on the current state $S_t$, chooses an action and appends a transition ($S_t, A_t, R_{t+1}, \gamma_{t+1}, S_{t+1}$) to a replay memory buffer that stores transitions. The parameters of the neural network are then optimized using stochastic gradient descent to minimize the loss.

\begin{equation}
\label{eq6}
\begin{aligned}
& L(\theta)=\left[R_{t+1}+\gamma_{t+1} \max _{a^{\prime}} Q_{\bar{\theta}}\left(S_{t+1}, a^{\prime}\right)-Q_\theta\left(S_t, A_t\right)\right]^2. \\
\end{aligned}
\end{equation}

The time step, $t$, is randomly selected from the replay memory. The loss gradient is then back-propagated exclusively into the parameters ($\theta$) of the online network, which is responsible for action selection. On the other hand, the term $\bar{\theta}$ represents the parameters of a target network, which serves as a periodically updated copy of the online network and is not directly optimized. To optimize the network, we employ Adam optimizer \citep{kingma2014adam}. We sample mini-batches uniformly from the experience replay and utilize them for the optimization process. 

\subsection{DQN Improvements}\label{section2_2}
The DQN algorithm has marked a significant milestone; however, it is important to acknowledge several limitations that affect its performance \citep{hessel2018rainbow}:

(a) Q-values Overestimation: The classical DQN algorithm generally suffers from overestimation bias, which arises due to the maximization step in Eqs. \ref{eq5}, \ref{eq6}) in which the update rule involves selecting the action with the maximum Q-value for the next state \citep{arulkumaran2017brief}. However, this includes the Q-values estimated by the same network being updated, which can lead to overestimation of the action values. The overestimation bias can result in suboptimal policies, where the agent tends to be overly optimistic about the value of certain actions, and therefore has adverse effects on learning \citep{hessel2018rainbow}. This bias increases as the complexity of the environment increases.

(b)Slow convergence:
The convergence rate of the classical DQN algorithm is slow when dealing with complex environments. During the early stages of training, the learning process can be unstable and may include oscillations \citep{sewak2019deep}. This instability arises primarily because of the interaction between function approximation (the neural network) and the use of the same network for action selection and evaluation. Many extensions have been proposed in the literature to address these limitations \citep{wan2020robust}. We integrate two well-known variants of the DQN algorithm to address the above issues.

-Double DQN: Double DQN aims to address the overestimation bias that exists in the classical DQN by employing two separate neural networks; the target network and the online network. The target network is dedicated to estimating the target Q-values, whereas the online network is responsible for the action selection \citep{haarnoja2018composable}. The separation of action selection and value estimation helps to reduce the overestimation bias \citep{hasselt2010double, van2016deep}. To incorporate Double Q-learning, the loss is calculated as follows:
\begin{equation}
\label{eq8}
\begin{split}
L(\theta) &=[R_{t+1}+\gamma_{t+1} Q_{\bar{\theta}}\left(S_{t+1} \underset{a^{\prime}}{\operatorname{argmax}} Q_\theta\left(S_{t+1}, a^{\prime}\right)\right) \\ 
&\left.-Q_\theta\left(S_t, A_t\right)\right]^2 .
\end{split}
\end{equation}
-Dueling DQN: Dueling DQN utilized an innovative network architecture dedicated to value-based model-free RL \citep{wang2016dueling}. It addresses the limitations of the classical DQN by separating the state value function estimation and the action advantages \citep{wang2016dueling}. This separation allows the network to understand which states are valuable, independent of the action chosen. The action values derived from the improved network are as follows:
\begin{equation}
\label{eq9}
\begin{split}
Q(s, a ; \xi, \eta, \psi) &= \mathrm{V}(s ; \xi, \eta) \\
&+\left(A(s, a ; \xi, \psi)-\frac{1}{N_{\text {act }}} \sum_{a^{\prime}} A\left(s, a^{\prime} ; \xi, \psi\right)\right).
\end{split}
\end{equation}
After all, Dueling DQN can generalize better across actions, leading to more efficient decision-making \citep{tavakoli2018action}. 
\subsection{DQN Algorithm}
\begin{figure*}
	\centering
	\includegraphics[trim=1cm 5cm 1cm 2cm, clip, scale=.4]{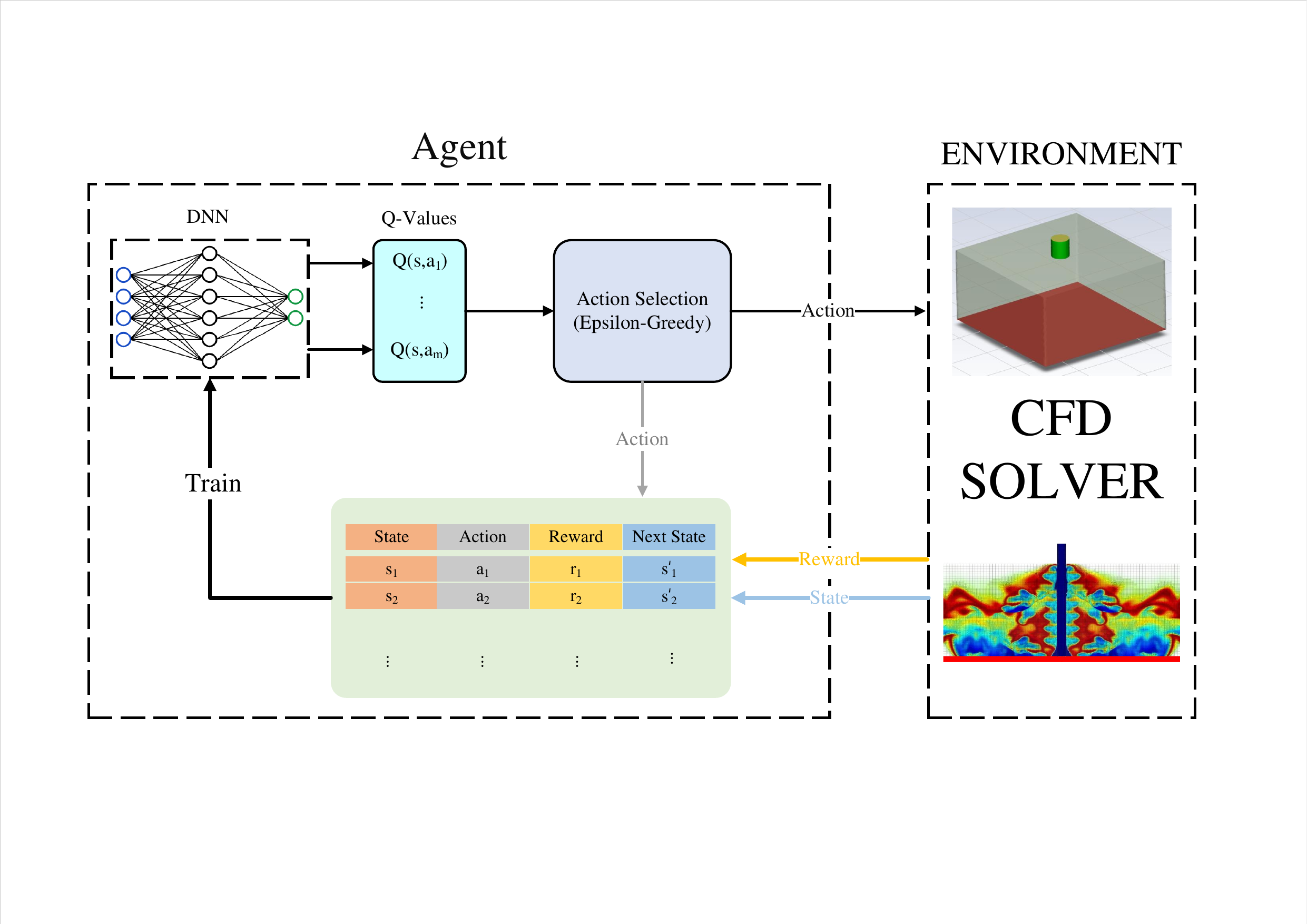}
	\caption{General overview of the DRL-CFD framework}
	\label{FIG_flowchart}
\end{figure*}

The real-time interaction between the CFD solver and the DRL agent is depicted in Fig. \ref{FIG_flowchart} and Algorithm 1. This framework can be divided into two main parts. In the first part, the state $S_t$ is formed by collecting data from the environment at time $t$. Subsequently, the agent chooses an action $\epsilon$-greedy, which is used by the CFD solver to solve the flow equations (Eq. \ref{eq1}) to the next state $S_t$. This interaction continues within each episode until termination.

Then, as for the second part, the sampled state $S_t$, along with the action $a_t$, reward $r_t$, and next state $S_{t+1}$, are stored in a replay buffer. The agent, then, extracts transitions of a specified size (mini-batches) and calculates the loss and updates the network parameters. These two parts run in parallel forming the deep Q-network. 

\section{Impingement Jet}

\subsection{Physics of the Flow}
\begin{figure}
	\centering
	\includegraphics[scale=.22]{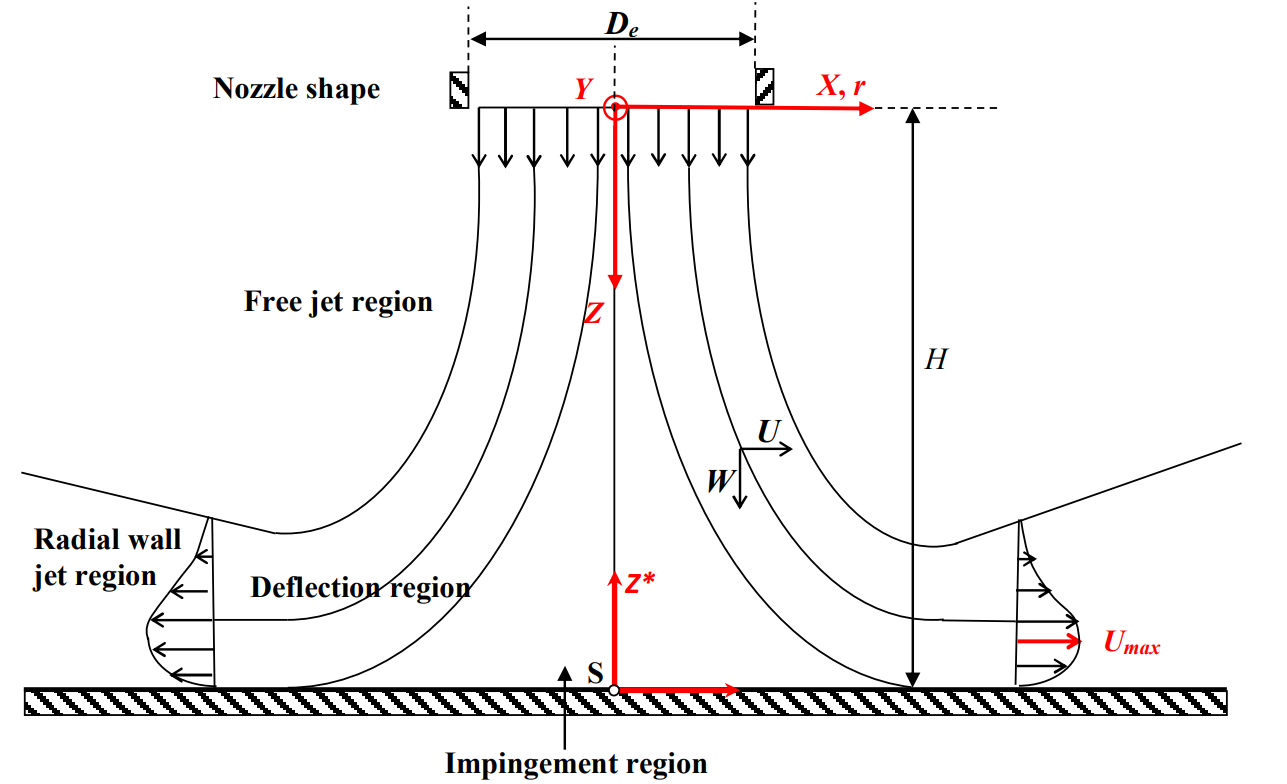}
	\caption{Formation of different regions in the context of the impinging jet on a flat plate \citep{sodjavi2015impinging} }
	\label{FIG_flowPhysics}
\end{figure}

Before discussing the DRL-CFD results, one needs to understand the physics of the flow in terms of impingement cooling. According to Fig. \ref{FIG_flowPhysics}, the flow field of an impinging jet consists of three regions; the free jet, the impingement region and the wall jet region. These regions are highly sensitive to the effective flow and geometrical parameters including Reynolds number, nozzle-to-plate distance, nozzle section shape and surface geometry. Change in any of the mentioned parameters results in changes in the boundaries between three regions, which subsequently also affects the heat transfer from the impinging surface. Among all the available parameters, the jet velocity can be considered as the most potential variable in active flow control of the impinging jet. A change in the jet velocity and frequency causes a change in the vortex shedding frequency in the flow field which eventually changes the formation and thickness of the boundary layer in the region of the wall jet. 

Moreover, in practical cases, it is difficult to find a general correlation that can provide an explicit relationship between the heat transfer coefficient or surface temperature in terms of jet velocity. Hence, the use of advanced machine learning techniques such as deep reinforcement learning is necessary for such systems.

\subsection{Setup Description}

\begin{figure*}
	\centering
	\includegraphics[scale=.1]{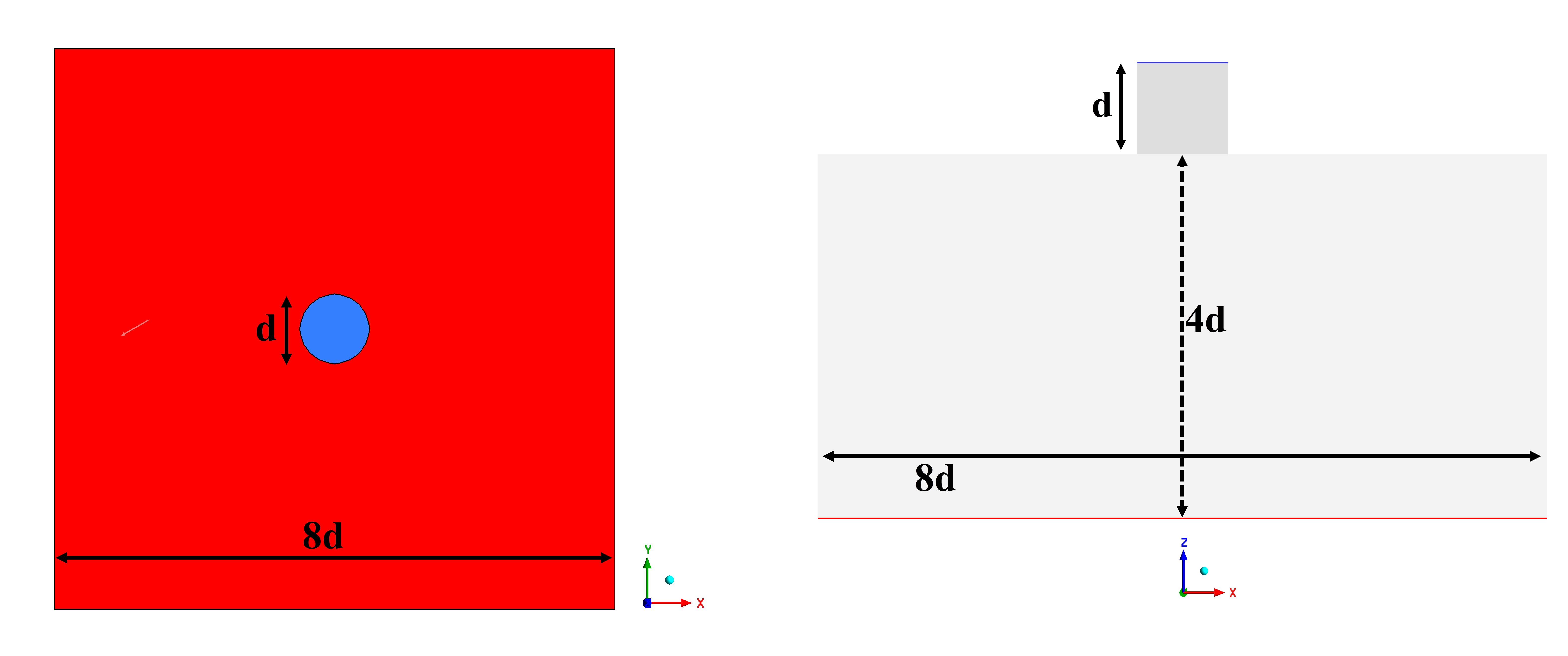}
	\caption{Schematic representation of the computational domain along with the dimension data of the jet and hot plate}
	\label{FIG_geom}
\end{figure*}

We focus on the forced convection of a 3-dimensional hot plate prone to a cooling jet. Fig. \ref{FIG_geom} shows the geometrical representation of the model. This model represents a simplified version of numerous complex cooling systems such as gas turbine cooling, electronics cooling, steel industry and aerospace applications \citep{bunker2014impingement, wang2020review, jena2022comparative}. The plate (shown in red) has a square shape with a length of 8d and is concentered with the impinging jet. It is subject to constant heat flux, $q ''$. Jet is placed to the distance of $H/d = 4$. Zero gradient Neumann boundary condition along with a constant pressure is considered for the outer sides of the domain. To decrease computational costs, we conduct calculations solely on one-quarter of the domain utilizing symmetry boundary conditions. The inlet of the jet (shown in blue) incorporates variable velocity magnitude ranging from $0.1V_\infty \sim V_\infty$ which is controlled by the DRL agent, while a constant temperature of $T_\infty$ enters the domain.

\begin{figure}
	\centering
	\includegraphics[scale=.06]{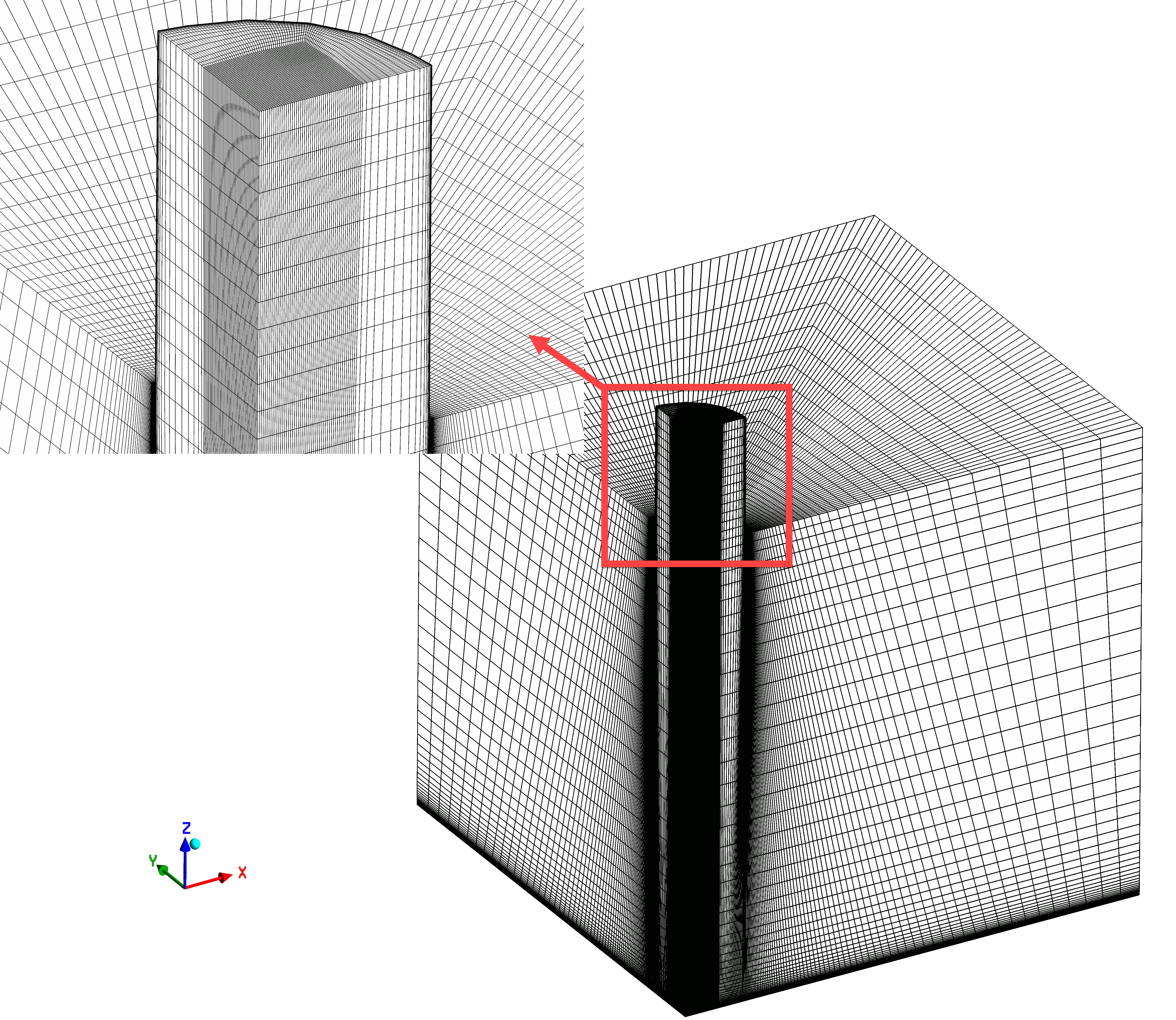}
	\caption{Representation of the structured mesh for the domain}
	\label{FIG_mesh}
\end{figure}

At the initial time, zero velocity and pressure were considered, while the isothermal temperature, $T_\infty$ is applied to the domain. For the sake of precision, a fully structured grid consisting of $8\times 10^5$ hexahedral elements has been utilized, as shown in Fig. \ref{FIG_mesh}.

The geometrical and flow parameters are inserted in Table \ref{table1}, along with fluid properties. Consequently, as for the variable velocity of the jet, the Reynolds number range will be:

\begin{equation}
\operatorname{Re}=\frac{\rho V d}{\mu}=170 \sim 1700
\end{equation}

which guarantees a laminar regime in the domain.

\begin{table*}[]
\caption{Parameters used in the 3-D thermal control}
\label{table1}
\begin{tabular}{cccccccc}
d     & $V_\infty$ &$T_\infty$   & $T_d$  & $\rho$                     & $\mu$       & $k$     & $C_p$   \\[2ex] \hline
($m$)   & ($m/s$) & ($K$) & ($K$)  & ($kg/m^3$) & ($Pa.s$)     & ($W/mK$)  & ($J/kg.K$) \\[1ex] \hline
0.025 & 1   & 288    & 303   & 1.225                   & 1.789e-5 & 0.024 & 1006  
\end{tabular}
\end{table*}

\section{Results and Disccusion}
The effectiveness of DRL-based thermal control of forced convection on the hot surface is initially validated by conducting a comparative analysis between the classical DQN and the baseline (i.e., no control situation) for the 3-D model. Subsequently, the control algorithm is expanded to explore variants of the DQN and provide a comparative study on the suitable approach for thermal control problems.

\subsection{DRL Specifications}
As mentioned in Section \ref{section2_1}, the agent's promotion is dictated by a reward-oriented strategy, where an optimized reward function enhances convergence. The following reward function was defined based on the averaged surface temperature, $T_{\text {surf }}$, and desired temperature, $T_{\text {d }}$:

\begin{equation}
\label{eq10}
\begin{cases}+1 & \left|T_{\text {surf }}-T_d\right|<2, \\ 0.1-\left(\left|T_{\text {surf }}/T_d-1\right|\right) \times 0.1 & \left|T_{\text {surf }}-T_d\right|>2.\end{cases}
\end{equation}

The training process consists of 100 episodes, each lasting 100 seconds ($10^{4}$ timesteps). Following training, the trained agent is exported for performance testing for 100 seconds.

\subsection{Sensitivity Analysis}
The appropriate state selection plays a crucial role in determining the agent's actions. Particularly, when dealing with high-dimensional problems, the state representation should encapsulate relevant information about the system's dynamics to reinforce the agent's prediction. The main challenge is to find a balance between the state representation richness and its computational possibility. Including too many variables might lead to a high-dimensional state space, which endangers convergence during training, whilst, an insufficient state representation can lead to a suboptimal solution. 

As for the selection of state variables, pressure, velocity, and temperature are the candidates. Among these, temperature is considered as a definitive choice due to its direct effect on thermal control. To determine the most suitable choice between velocity and pressure, we examined their temporal/spatial gradient within the domain. It was observed that velocity exhibits higher gradients all over the domain due to the hydraulic boundary layer. Thus, velocity and temperature can be served as the primary state variables for thermal control systems. 

\begin{figure}
	\centering
	\includegraphics[scale=.22]{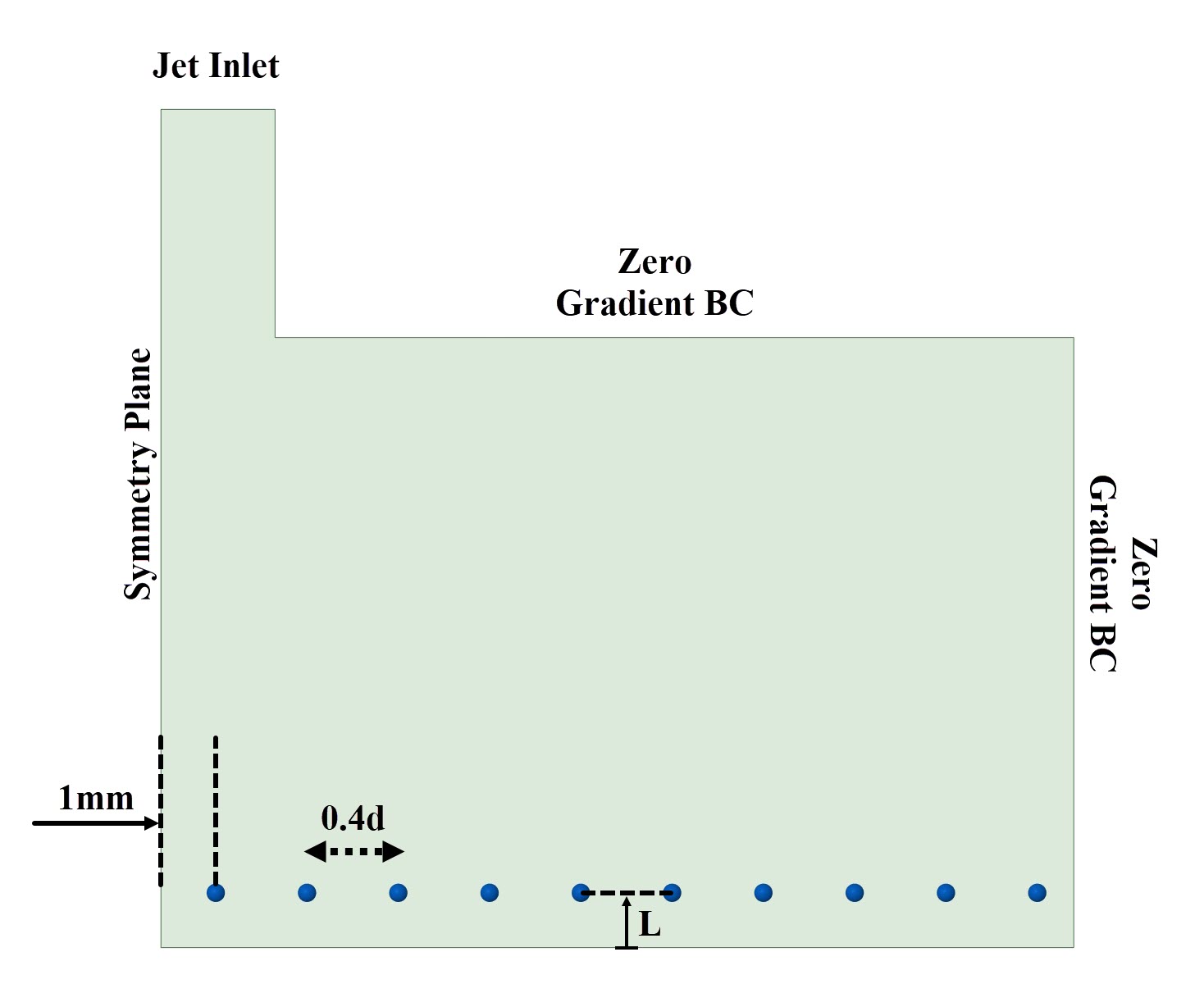}
	\caption{Probes location for measurement of states during DRL are shown with blue dots}
	\label{FIG_probes}
\end{figure}

The next complexity is to determine the suitable sensor location for measuring state variables. A sensitivity analysis has been performed by considering three distances between the probes and the surface, denoted as L in Fig. \ref{FIG_probes}: $1$ mm, $5$ mm, and $10$ mm. 

\begin{figure}
	\centering
	\includegraphics[scale=.3]{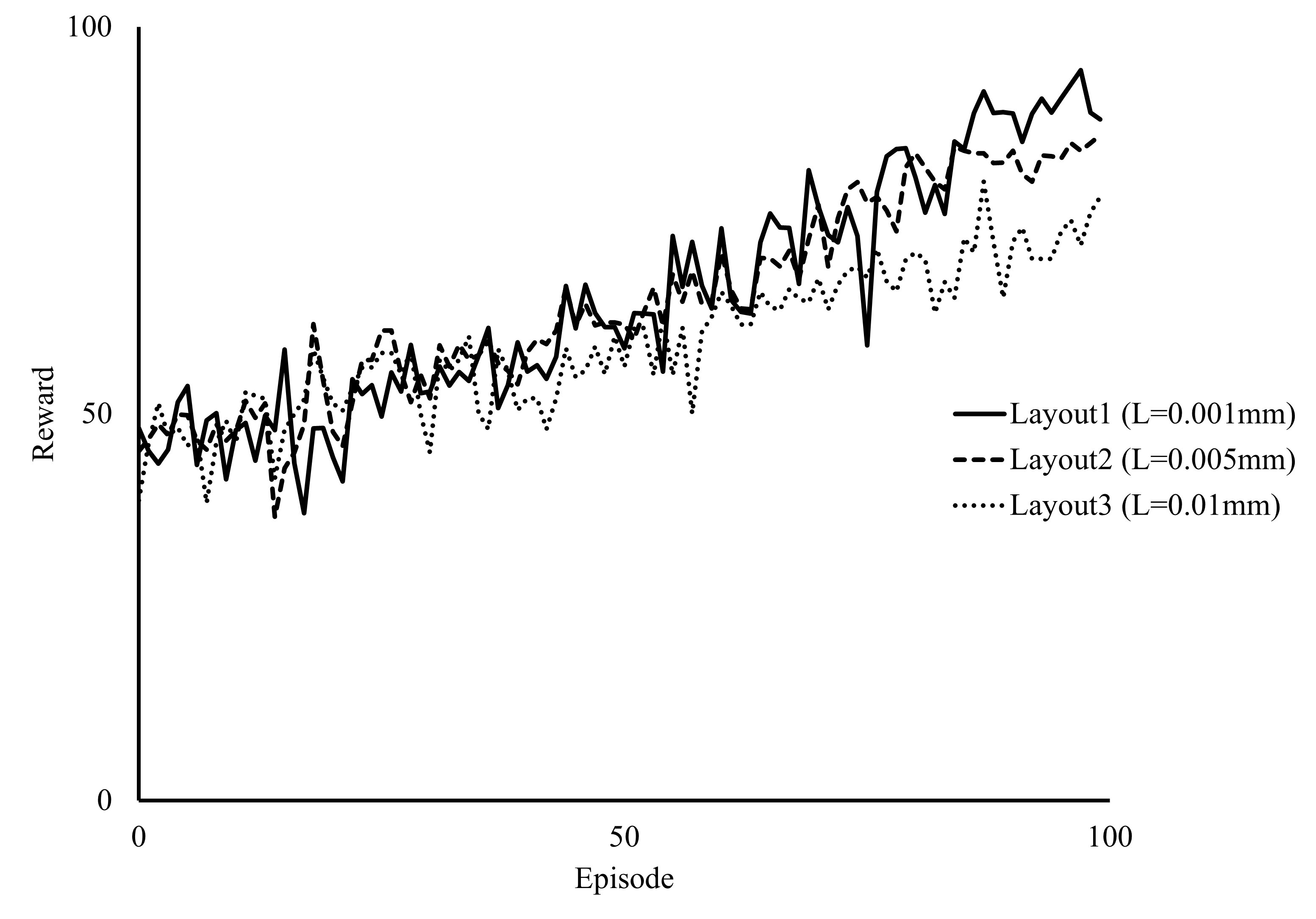}
	\caption{Comparison between total reward for different sensor layouts}
	\label{FIG_probesReward}
\end{figure}

Fig. \ref{FIG_probesReward} shows the change of total reward in terms of episode number for the mentioned layouts. As depicted, layouts $1$ and $2$ exhibit similar trends while the former demonstrated superior performance at the end. Layout1 corresponds to the placement of the sensors closest to the wall, i.e. $L=0.001$ mm.

\subsection{Classical DQN}
\begin{figure}
	\centering
	\includegraphics[scale=.32]{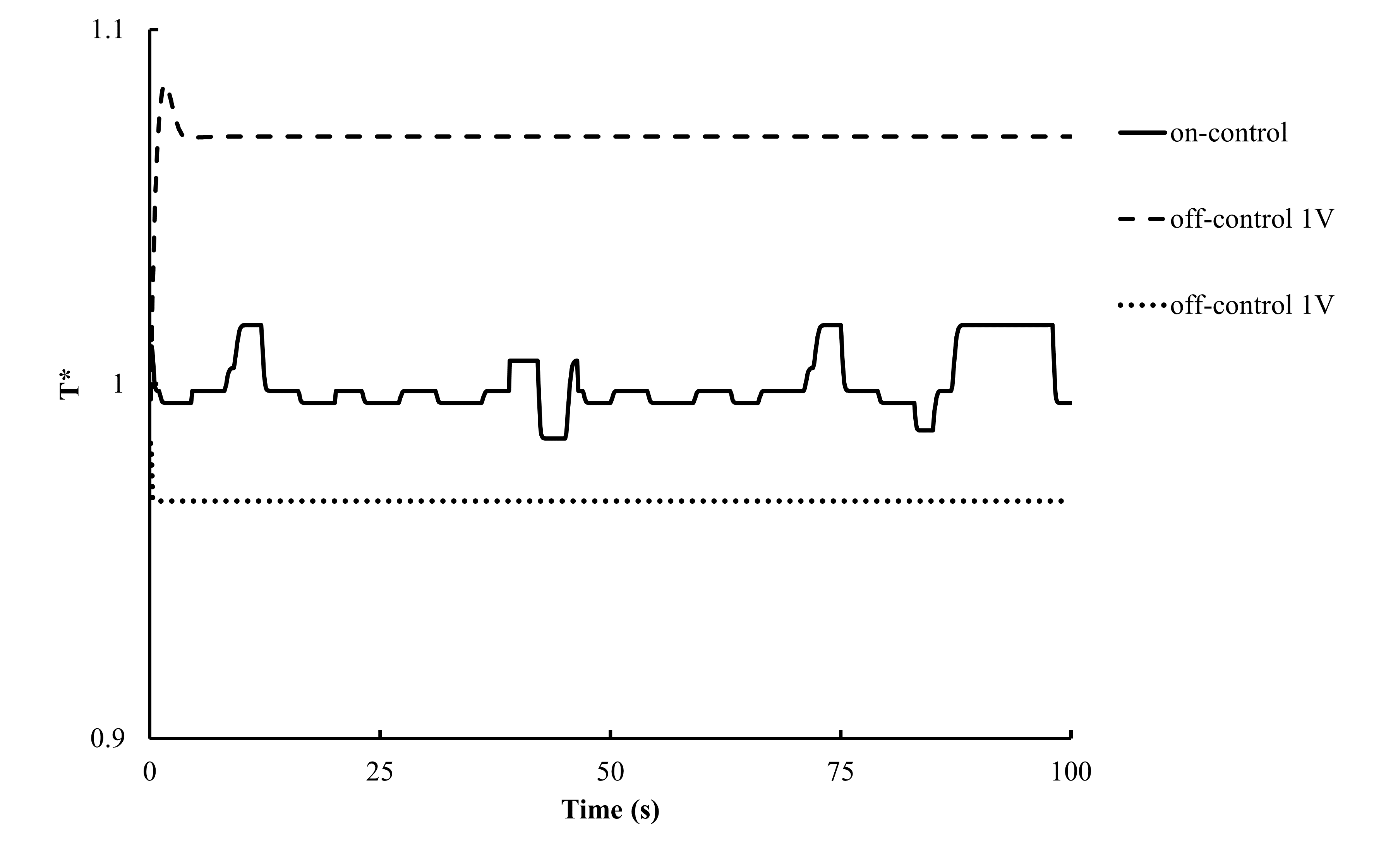}
	\caption{Temporal variation of dimensionless surface average temperature for off- and on-control (classic DQN) systems}
	\label{FIG_off_onControlPlot}
\end{figure}

As the first part of this research study, we aim to assess the applicability of Deep Reinforcement Learning for thermal control based on a custom environment defined by a CFD solver. Thus, we consider the agent, represented by the jet velocity, within a velocity range of $0.1V_\infty$ to $V_\infty$. Fig. \ref{FIG_off_onControlPlot} shows the comparison between the temporal changes of dimensionless average surface temperature ($T^*=T/T_d$) for the lower and higher velocity ranges without any control, along with the results obtained using a DQN-based velocity control. The uncontrolled strategy exhibits average temperatures of $1.07$ ($324$ K) and $0.97$ ($294$ K) for the lowest and highest velocities, respectively. In contrast, by applying the DRL-CFD method, we successfully maintained the dimensionless average surface temperature at an almost average of $1$ ($303$ K). These results indicate the capability of DQN to achieve effective thermal control systems.

\begin{figure*}
	\centering
	\includegraphics[scale=.08]{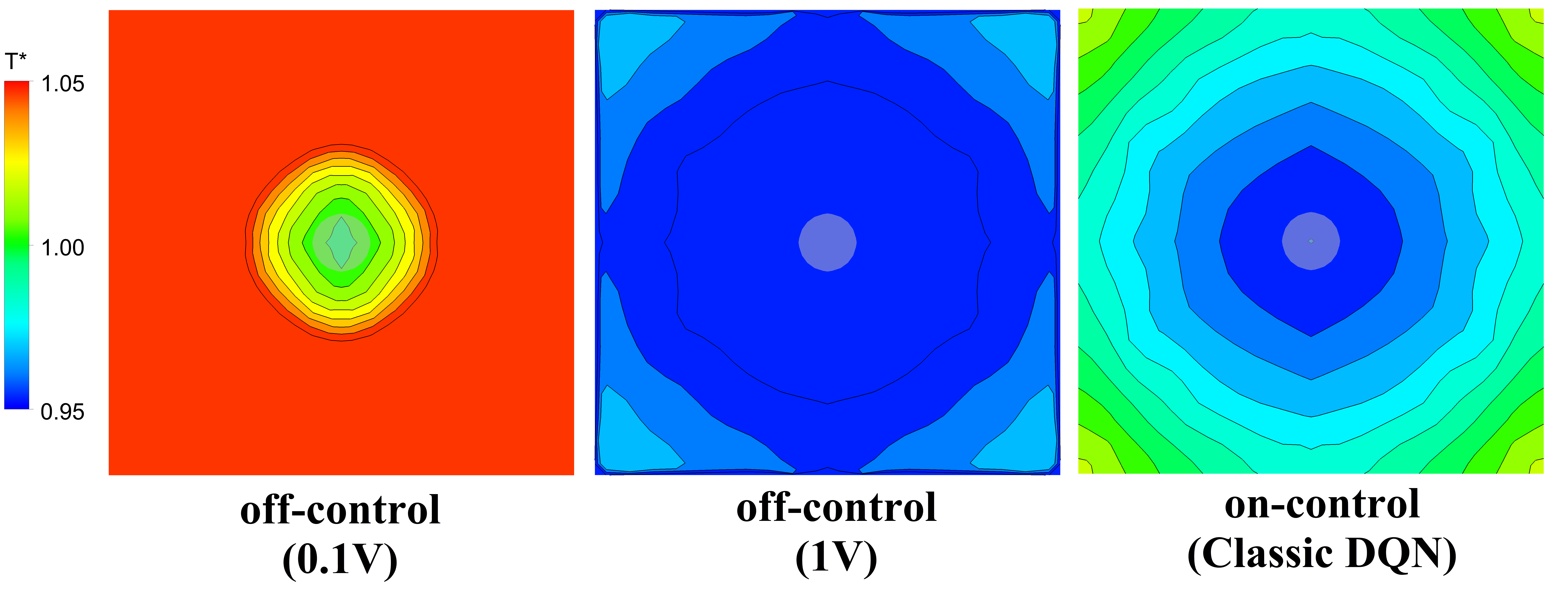}
	\caption{Temperature contour for off- and on-control (classic DQN) systems}
	\label{FIG_off_onControlContour}
\end{figure*}

The contour plot presented in Fig. \ref{FIG_off_onControlContour} illustrates the average temperature distribution in off- and on-control systems. In the case of the uncontrolled system with a low jet velocity, elevated temperature values are observed throughout the entire surface, except for the central region where jet flow occurs directly. Conversely, the application of a high jet velocity leads to significantly lower temperatures across the surface.

For the controlled system, the contour plot is taken by statistically averaging temperature from $0$ $\sim$ $100$ seconds. Remarkably, the controlled system consistently maintains an average dimensionless temperature of $1$ across the entire surface, with a deviation of less than $2$\%. These findings highlight the effective temperature regulation achieved by the thermal control system, ensuring stability and uniformity throughout the system.

\subsection{Episode Number}
\begin{figure}
    \centering

    \begin{subfigure}[b]{0.45\textwidth}
        \centering
        \includegraphics[width=\textwidth]{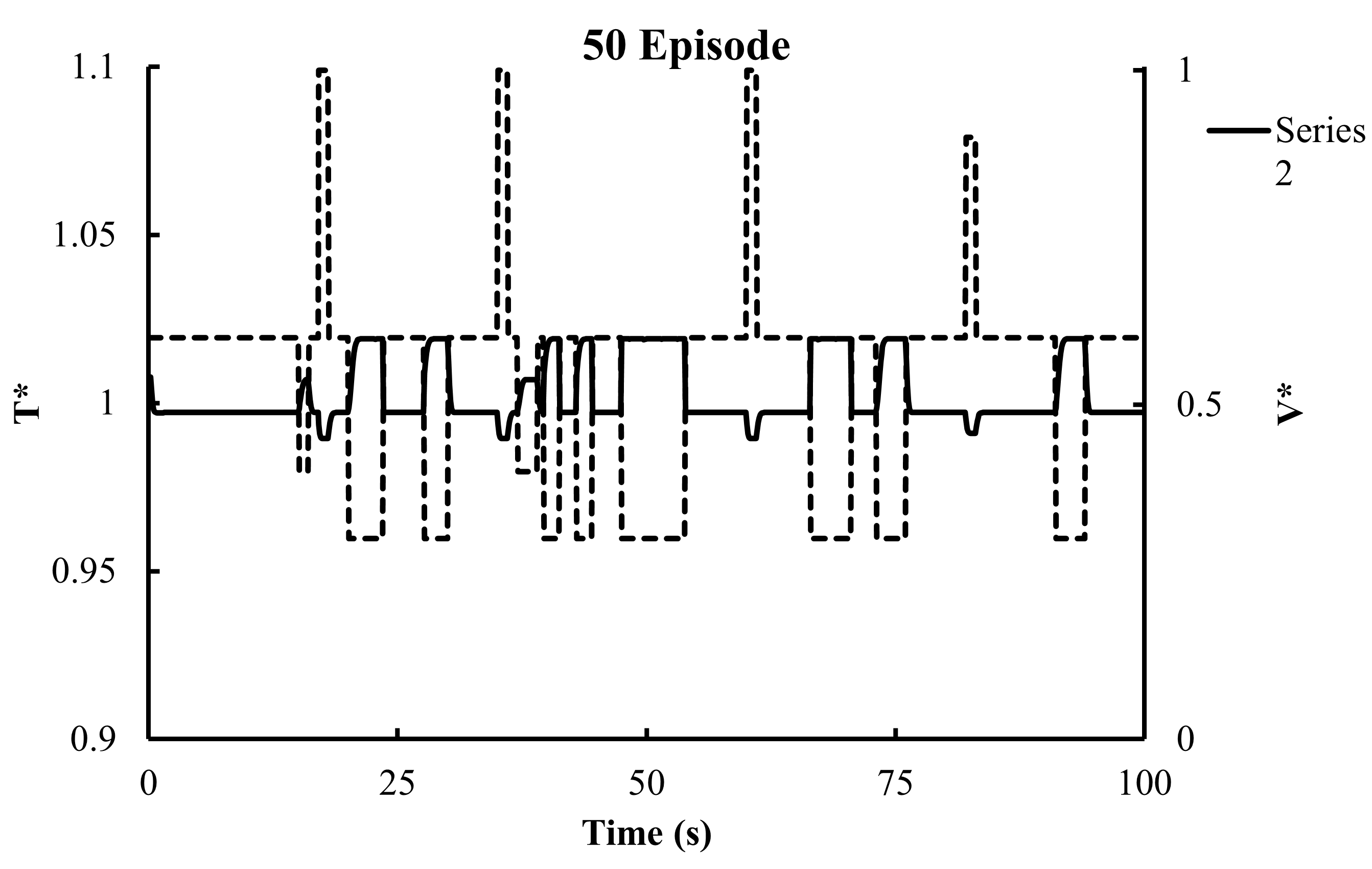}
    \end{subfigure}
    \hfill
    \begin{subfigure}[b]{0.45\textwidth}
        \centering
        \includegraphics[width=\textwidth]{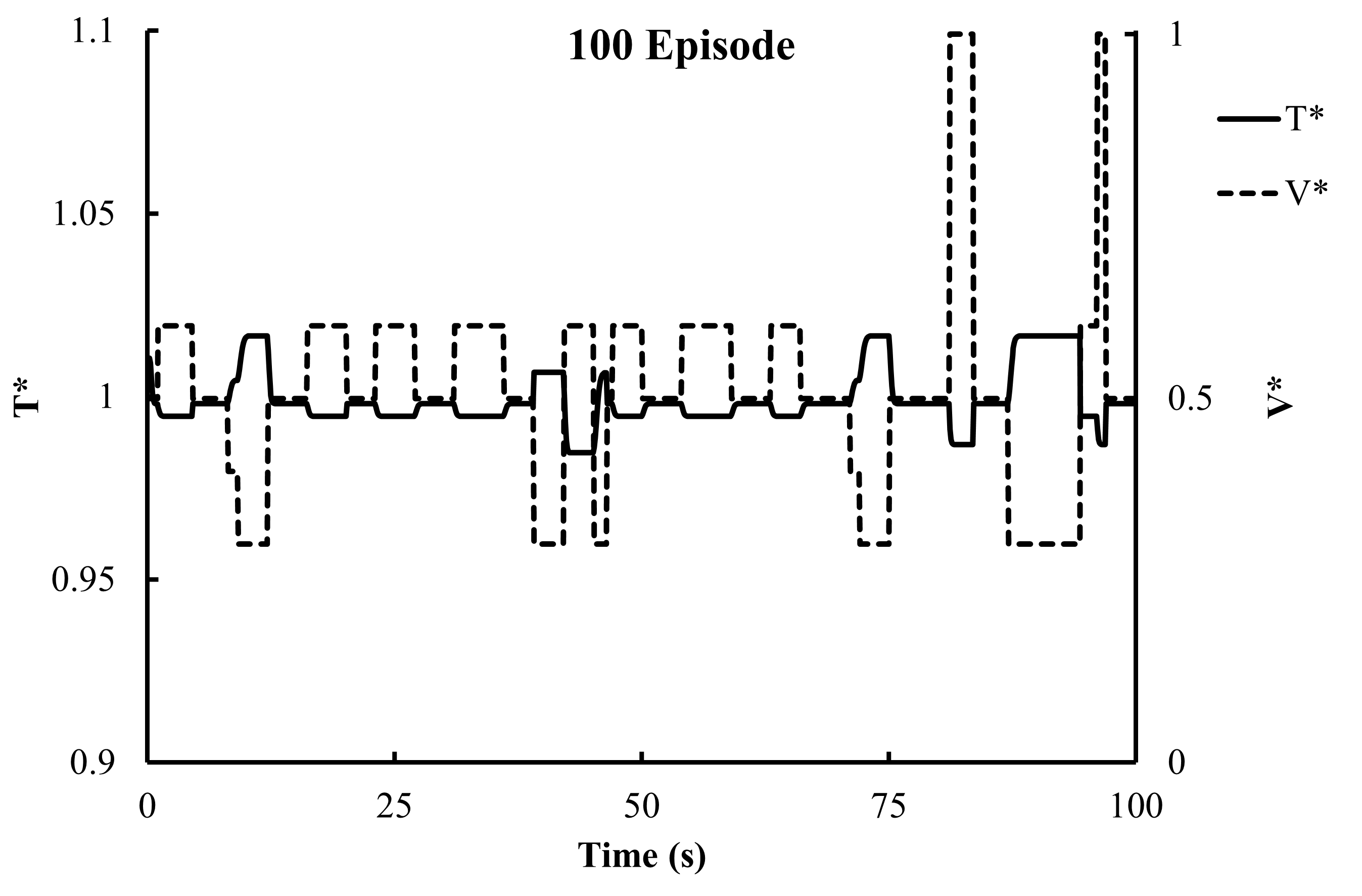}
    \end{subfigure}
    \hfill
    \begin{subfigure}[b]{0.45\textwidth}
        \centering
        \includegraphics[width=\textwidth]{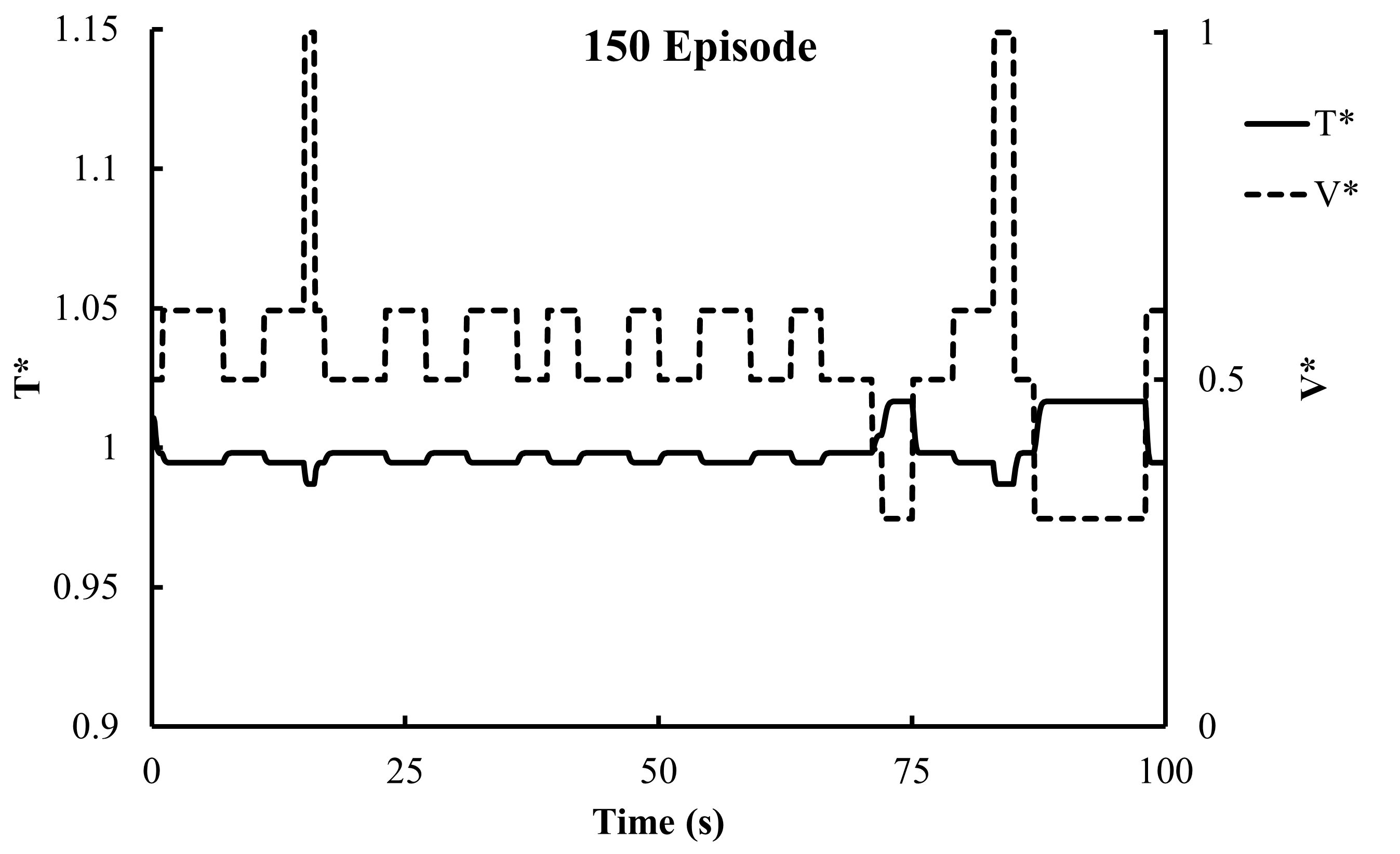}
    \end{subfigure}

    \caption{Time history of surface temperature and jet velocity for different episode numbers}
    \label{FIG_off_nEpisodeIndep}
\end{figure}

We conducted three training runs with varying episode numbers: $50$, $100$, and $150$. Fig. \ref{FIG_off_nEpisodeIndep} presents the surface average temperature based on the on-control agents trained by different episode numbers. The agent trained with $50$ episodes shows higher levels of oscillatory behaviour when utilized for the on-control purpose. However, those trained by $100$ and $150$ episodes demonstrated lower oscillation, indicating a more desirable behaviour for a controller.

\begin{figure*}
	\centering
	\includegraphics[scale=.08]{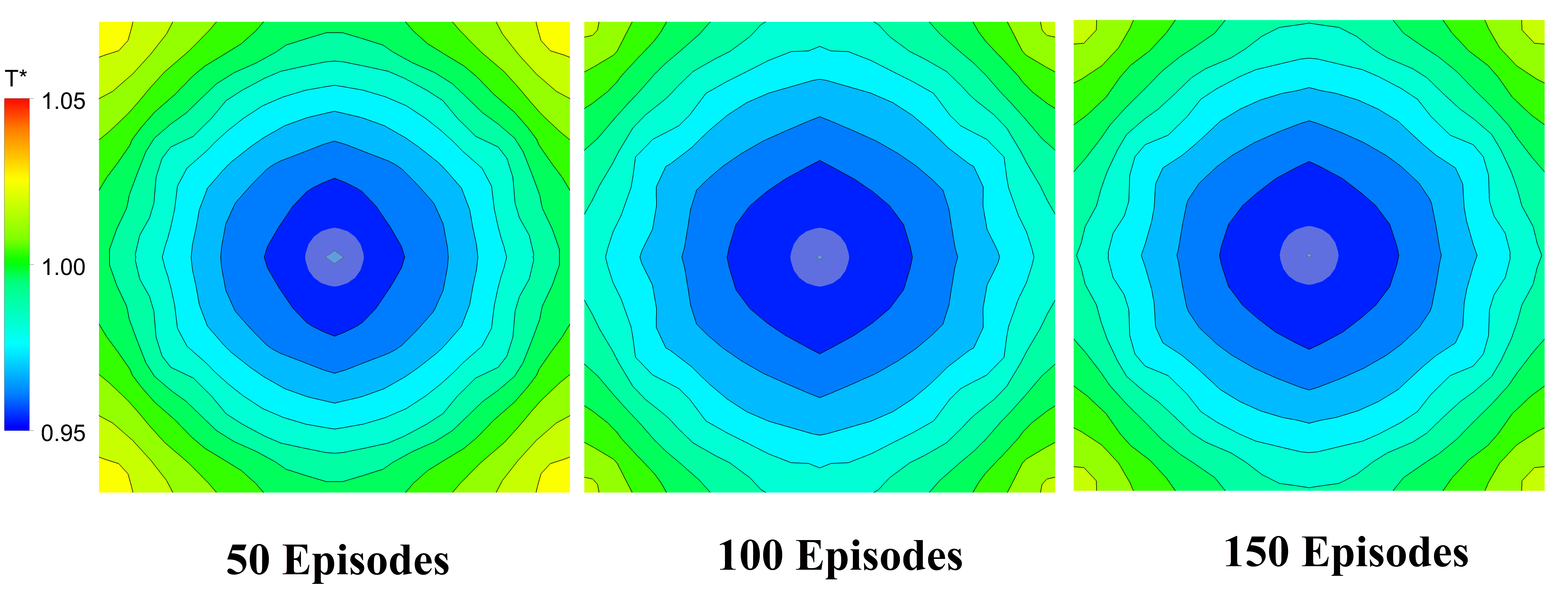}
	\caption{Temperature contour for different episode numbers}
	\label{FIG_nEpisodeUndepContour}
\end{figure*}

Fig. \ref{FIG_nEpisodeUndepContour}, displays the time-averaged temperature contour on the hot plate for three episode numbers. The temperature contour observed in episodes 100 and 150 exhibit similarities, with cooler temperatures at the center of the surface ($<T_d$) and higher temperatures ($>T_d$) at the corners. However, episode 50 demonstrates comparatively lower performance as it incorporates higher temperature with respect to other episodes.

\subsection{DQN Variants}
In this part, we will conduct a comparative analysis of different variants of DQN. As discussed in Section \ref{section2_2}, our focus will be on investigating the effectiveness of classical DQN, Double DQN, and Duel DQN. 

The main idea behind Double Q-learning is to lessen the overestimation bias that can occur with classical Q-learning algorithms. To address this issue, two subvariants of this method have been developed in terms of target network updating; hard update and soft update. In the hard Double DQN, the update of the target network occurs at specific intervals (i.e. after a fixed number of steps). Whereas, soft Double DQN introduces a more gradual way of updating the target network, in which instead of copying the entire parameters from the \textit{online} network, the soft update performs a weighted average between the online network and the target network as follows \citep{wei2018multiagent}:

\begin{equation}
\label{eq11}
\theta_\text{{target}} = \tau \theta_\text{{online}} + (1-\tau) \theta_\text{{target}},
\end{equation}

where $\theta_\text{{target}}$ and $\theta_\text{{online}}$ refer to the parameters of the target and online networks, respectively. $\tau$ is a small value close to zero (here we take 0.001) that controls the extent of the update.

\begin{figure}
	\centering
	\includegraphics[scale=.45]{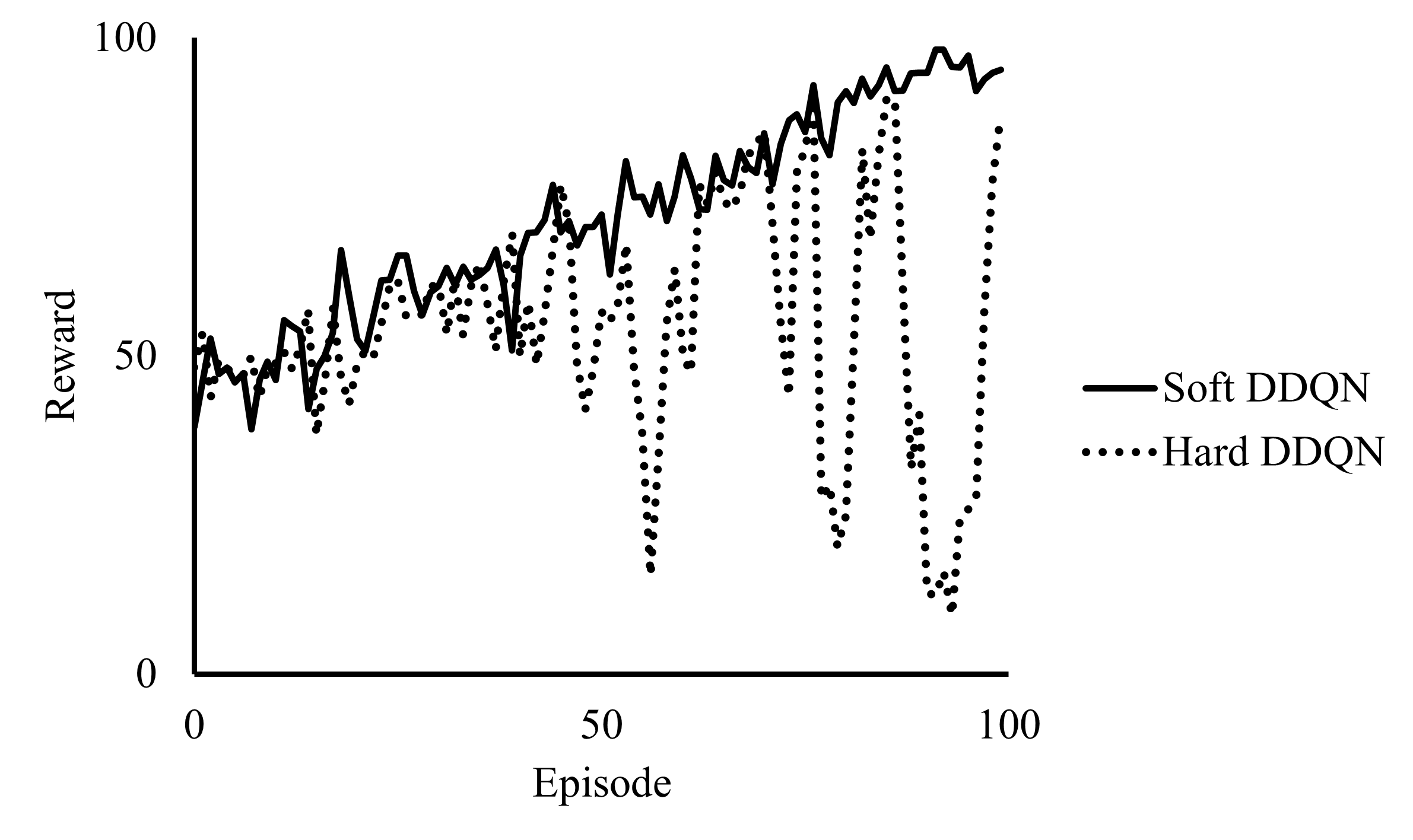}
	\caption{Evolution of the total reward value for soft and hard Double DQN in terms of the episode number}
	\label{FIG_variantsDDQN}
\end{figure}

Fig. \ref{FIG_variantsDDQN} illustrates the evolution of total rewards for both soft and hard Double DQN. The hard DQN exhibits high-frequency oscillations with high amplitude, indicating its failure in the thermal control task. However, the soft Double DQN demonstrates a substantial increase in total reward across episodes, reaching a remarkable value of 95 at the end. It is worth mentioning that a total reward of 95 implies the achievement of the desired temperature in more than 95\% of time instances during the final episode (based on the reward definition in Eq. \ref{eq10}), demonstrating its excellent performance.

\begin{figure*}
	\centering
	\includegraphics[scale=.2]{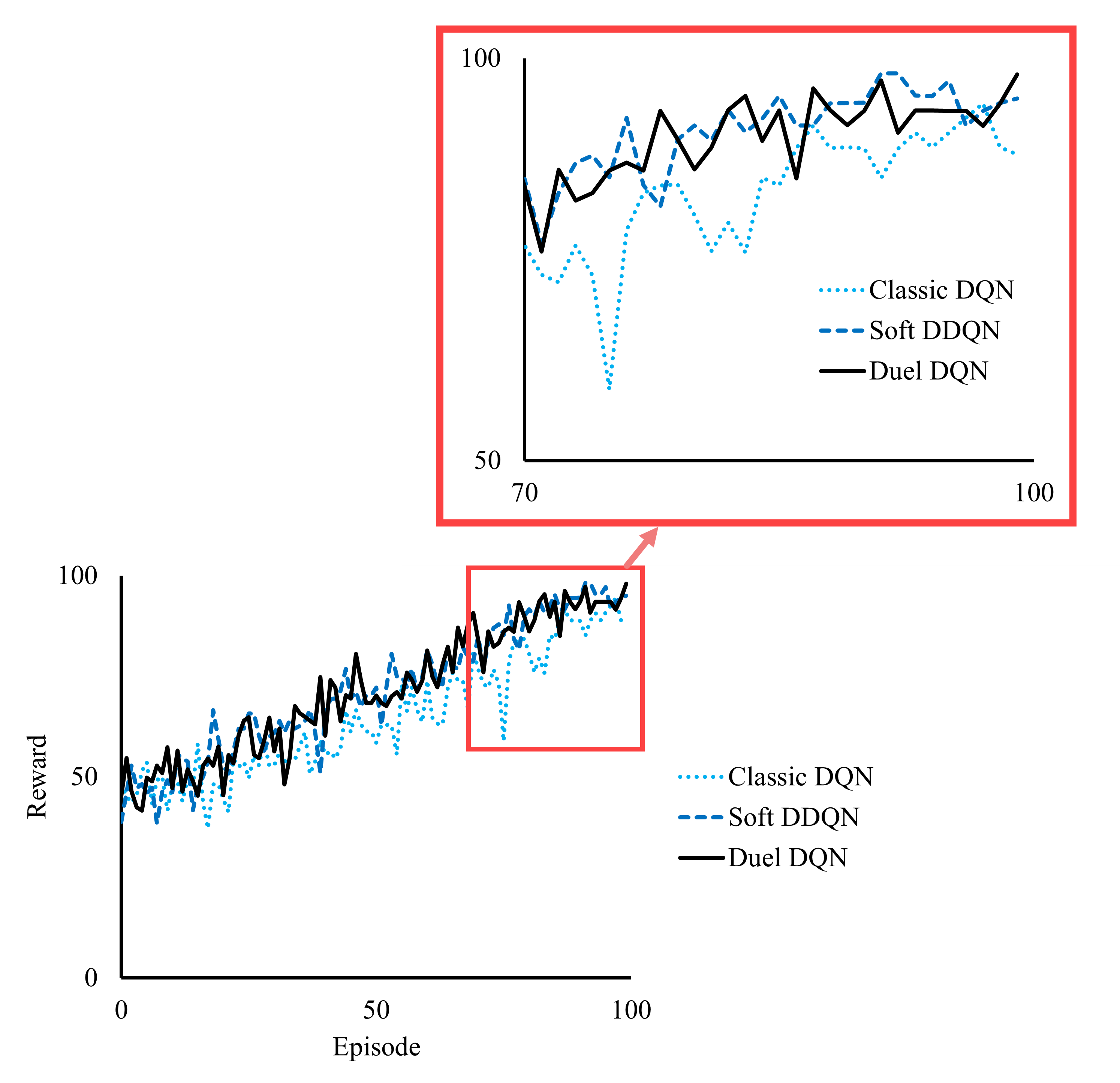}
	\caption{Evolution of the total reward value for different variants of DQN in terms of the episode number}
	\label{FIG_variantsReward}
\end{figure*}

Fig. \ref{FIG_variantsReward} provides a comparison of the reward evolution for DQN variants. Overall, the different variants perform similarly in terms of thermal control. However, the classical DQN shows a lower reward throughout the episodes and demonstrates more oscillatory behaviour with respect to Double and Duel DQNs. Upon a closer examination of the last 30 episodes for the three variants, Duel DQN shows the best final reward, followed by the Soft Double DQN and then Classical DQN.

\begin{figure*}
	\centering
	\includegraphics[scale=.11]{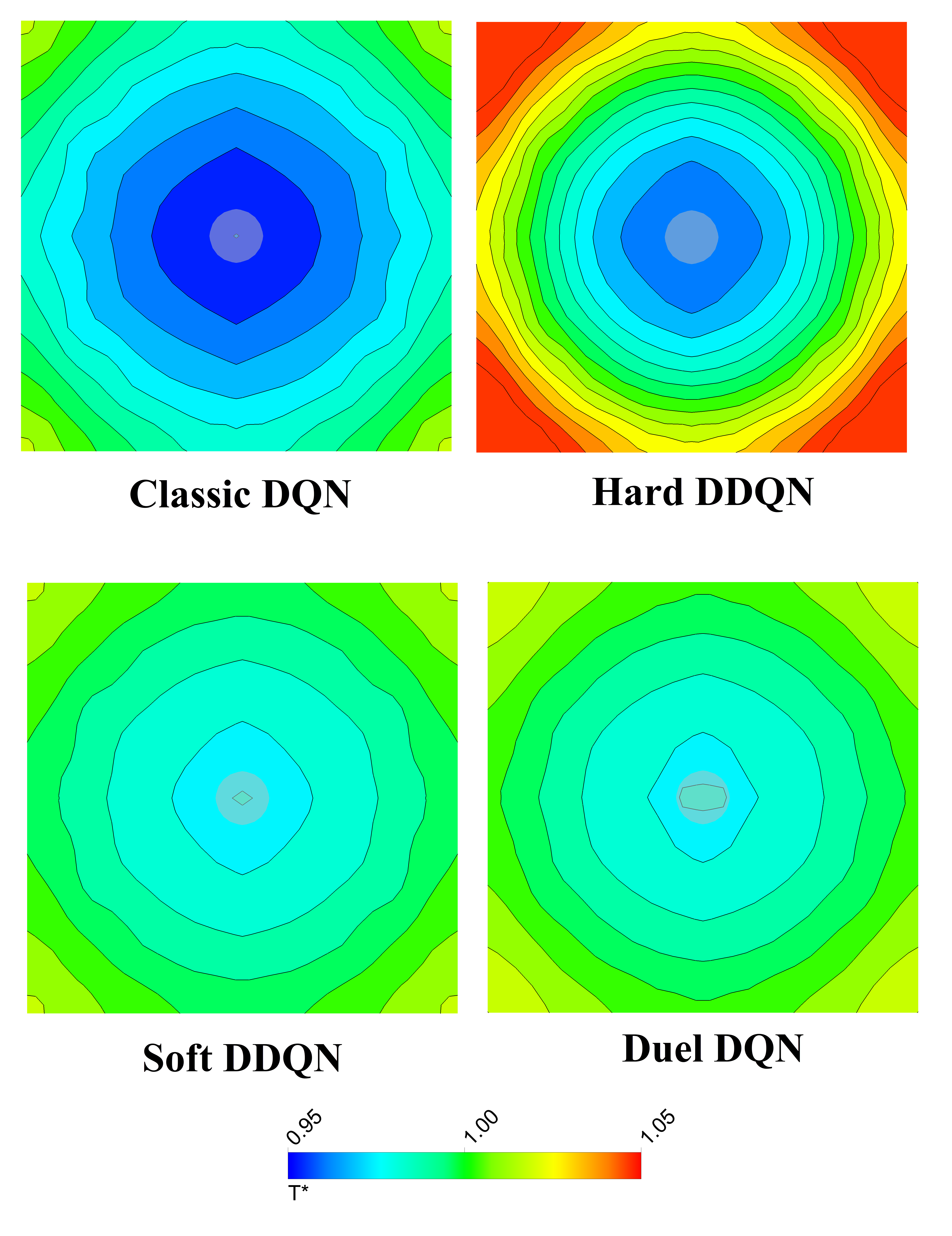}
	\caption{Temperature contour for different variants of DQN}
	\label{FIG_ContComp_variants}
\end{figure*}

By analysing the time-averaged surface temperature for all the agents, Fig. \ref{FIG_ContComp_variants} reveals that the classical DQN exhibits higher temperature gradients on the surface. This behaviour corresponds to the fact that the agent (jet velocity) experiences higher changes over time. On the other hand, both Soft Double DQN and Duel DQN achieve a more uniform surface temperature close to the desired value of $T^*=1$ (303 K). Moreover, the contour of Hard DQN illustrates the ineffectiveness of this method in thermal control systems.

\begin{figure}
    \centering

    \begin{subfigure}[b]{0.45\textwidth}
        \centering
        \includegraphics[width=\textwidth]{figs/EIndep_E100.png}
    \end{subfigure}
    \hfill
    \begin{subfigure}[b]{0.45\textwidth}
        \centering
        \includegraphics[width=\textwidth]{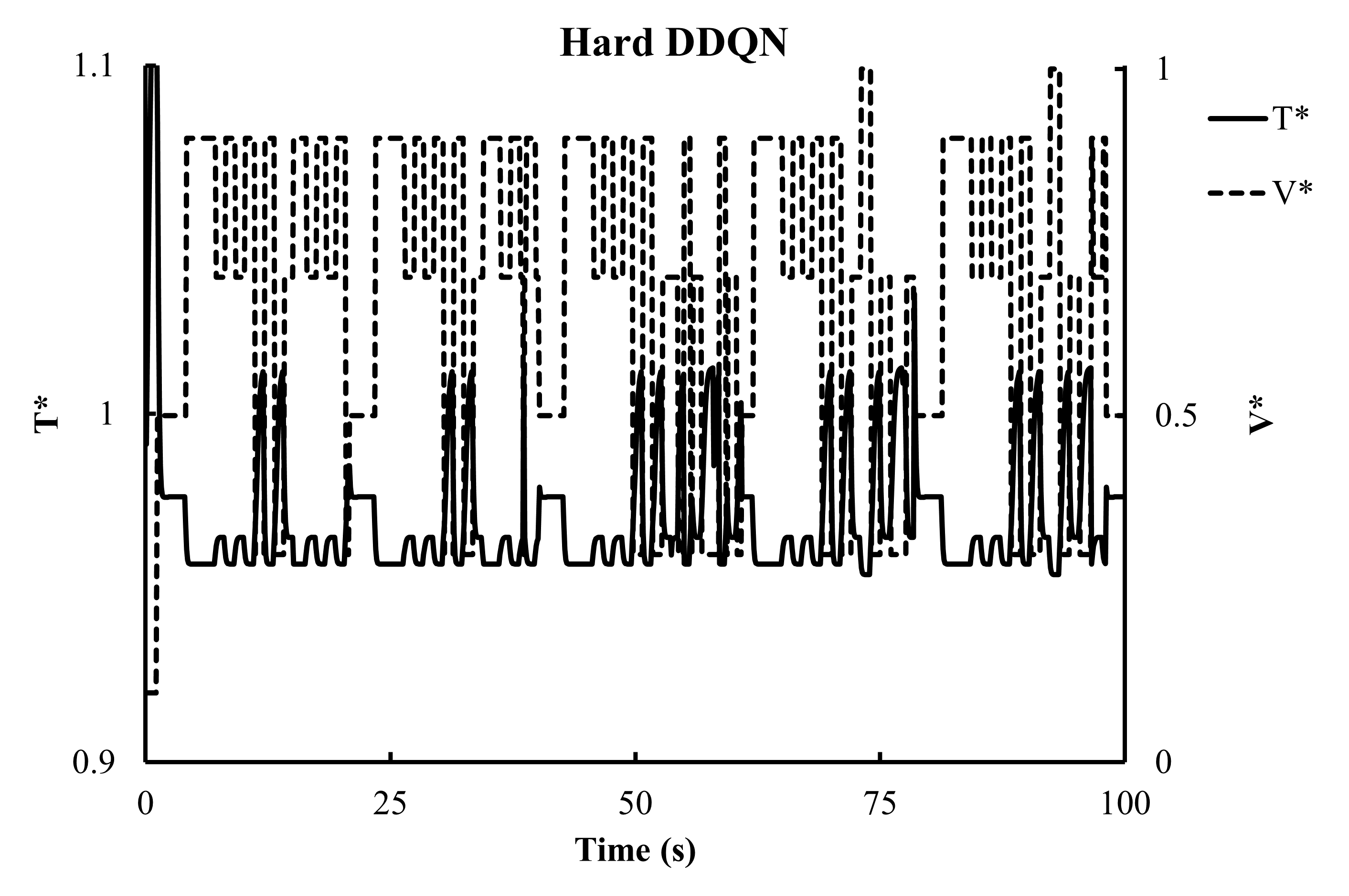}
    \end{subfigure}
    \hfill
    \begin{subfigure}[b]{0.45\textwidth}
        \centering
        \includegraphics[width=\textwidth]{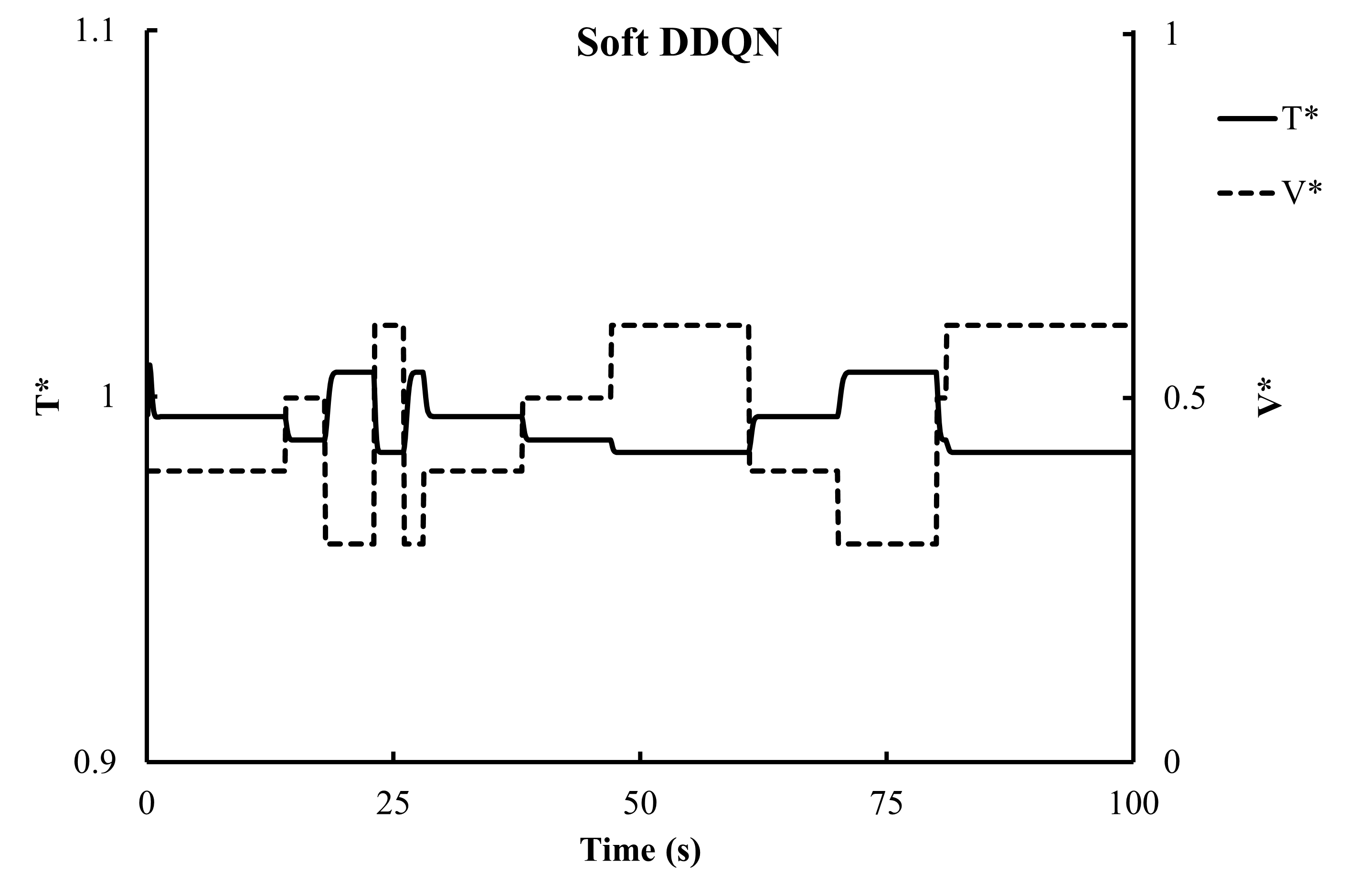}
    \end{subfigure}
    \hfill
    \begin{subfigure}[b]{0.45\textwidth}
        \centering
        \includegraphics[width=\textwidth]{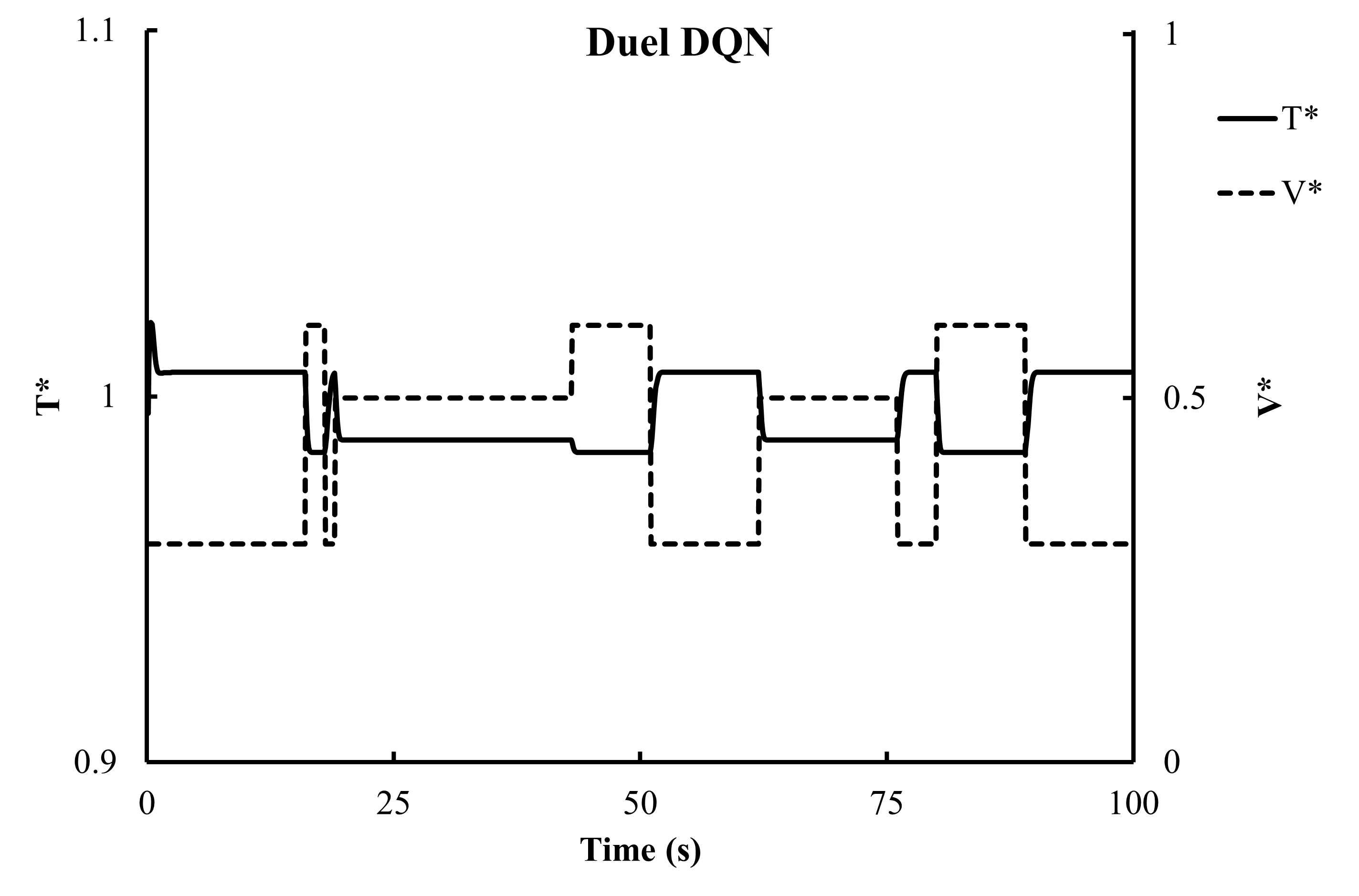}
    \end{subfigure}

    \caption{Time history of surface temperature and jet velocity for on-control simulation of the DQN variants}
    \label{FIG_variantVandT}
\end{figure}

Fig. \ref{FIG_variantVandT} presents the changes in temperature and velocity during the on-control simulation for DQN variants. The Soft Double DQN and Duel DQN exhibit minimal velocity changes, resulting in a more uniform temperature distribution close to the desired setpoint, as predicted. However, Hard DQN's significant changes in actions (jet velocity) lead to poor temperature control, with a considerable discrepancy from the desired temperature.

The Soft Double DQN and Dueling DQN methods prove to be effective for jet velocity control, as they produce stable and close-to-desired temperatures on the surface. The Classical DQN and Hard DQN variants, on the other hand, exhibit limitations in achieving optimal thermal control performance. These findings highlight the importance of selecting appropriate DQN variants for thermal control problems to achieve the desired temperature efficiently and reliably.

\section{Conclusion}
This research study focused on the application of the Deep Q-Network (DQN) method for a thermal control problem consisting of a hot plate prone to a cooling jet agent with the aid of a custom CFD environment. Two subvariants of Double DQN, namely hard update and soft update, were explored to address the overestimation bias that can occur with classical Q-learning algorithms. The results clearly demonstrate that the soft Double DQN outperforms the hard Double DQN in achieving effective thermal control on the hot plate.

Comparison with classical DQN and its advanced variants revealed that both soft Double DQN and Dueling DQN methods performed well in achieving stable and close-to-desired temperature control. The classical and hard DQN variants, however, showed limitations in achieving optimal thermal control, with the classical DQN displaying oscillatory behaviour.

By analyzing the time-averaged surface temperature, it was observed that the classical DQN exhibited higher temperature gradients along the surface, corresponding to higher jet velocity changes through the control cycle. In contrast, both soft Double DQN and Dueling DQN achieved a more uniform surface temperature close to the desired value.

Although all the variants except for the hard Double DQN demonstrated an acceptable performance for the thermal control task, the soft Double DQN and Dueling DQN methods proved to be more effective for jet velocity control, producing stable and close-to-desired temperatures. Future research could explore other advanced deep reinforcement learning algorithms to further enhance thermal control systems.

\bibliographystyle{unsrtnat}

\bibliography{refs}

\begin{thebibliography}{64}
\providecommand{\natexlab}[1]{#1}
\providecommand{\url}[1]{\texttt{#1}}
\expandafter\ifx\csname urlstyle\endcsname\relax
  \providecommand{\doi}[1]{doi: #1}\else
  \providecommand{\doi}{doi: \begingroup \urlstyle{rm}\Url}\fi

\bibitem[Bergman et~al.(2011)Bergman, Lavine, Incropera, and DeWitt]{bergman2011introduction}
Theodore~L Bergman, Adrienne~S Lavine, Frank~P Incropera, and David~P DeWitt.
\newblock \emph{Introduction to heat transfer}.
\newblock John Wiley \& Sons, 2011.

\bibitem[Childs et~al.(1999)Childs, Greenwood, and Long]{childs1999heat}
PRN Childs, JR~Greenwood, and CA~Long.
\newblock Heat flux measurement techniques.
\newblock \emph{Proceedings of the Institution of Mechanical Engineers, Part C: Journal of Mechanical Engineering Science}, 213\penalty0 (7):\penalty0 655--677, 1999.

\bibitem[Bejan(2013)]{bejan2013convection}
Adrian Bejan.
\newblock \emph{Convection heat transfer}.
\newblock John wiley \& sons, 2013.

\bibitem[Giwa et~al.(2021)Giwa, Sharifpur, Ahmadi, and Meyer]{giwa2021review}
SO~Giwa, M~Sharifpur, MH~Ahmadi, and JP~Meyer.
\newblock A review of magnetic field influence on natural convection heat transfer performance of nanofluids in square cavities.
\newblock \emph{Journal of Thermal Analysis and Calorimetry}, 145:\penalty0 2581--2623, 2021.

\bibitem[Xiong et~al.(2021)Xiong, Hamid, Iqbal, Irfan, and Khan]{xiong2021numerical}
Pei-Ying Xiong, Aamir Hamid, Kaleem Iqbal, M~Irfan, and Masood Khan.
\newblock Numerical simulation of mixed convection flow and heat transfer in the lid-driven triangular cavity with different obstacle configurations.
\newblock \emph{International Communications in Heat and Mass Transfer}, 123:\penalty0 105202, 2021.

\bibitem[Shahabadi et~al.(2021)Shahabadi, Mehryan, Ghalambaz, and Ismael]{shahabadi2021controlling}
Mohammad Shahabadi, SAM Mehryan, Mohammad Ghalambaz, and Muneer Ismael.
\newblock Controlling the natural convection of a non-newtonian fluid using a flexible fin.
\newblock \emph{Applied Mathematical Modelling}, 92:\penalty0 669--686, 2021.

\bibitem[Davalath and Bayazitoglu(1987)]{davalath1987forced}
J~Davalath and Y~Bayazitoglu.
\newblock Forced convection cooling across rectangular blocks.
\newblock 1987.

\bibitem[Al-Sarkhi and Abu-Nada(2005)]{al2005characteristics}
Abdelsalam Al-Sarkhi and E~Abu-Nada.
\newblock Characteristics of forced convection heat transfer in vertical internally finned tube.
\newblock \emph{International Communications in Heat and Mass Transfer}, 32\penalty0 (3-4):\penalty0 557--564, 2005.

\bibitem[{\"O}ztop et~al.(2006)]{oztop2006turbulence}
Hakan~F {\"O}ztop et~al.
\newblock Turbulence forced convection heat transfer over double forward facing step flow.
\newblock \emph{International communications in heat and mass transfer}, 33\penalty0 (4):\penalty0 508--517, 2006.

\bibitem[Kim et~al.(2022)Kim, Pokharel, and Yeom]{kim2022enhancing}
Kiyun Kim, Pravesh Pokharel, and Taiho Yeom.
\newblock Enhancing forced-convection heat transfer of a channel surface with synthetic jet impingements.
\newblock \emph{International Journal of Heat and Mass Transfer}, 190:\penalty0 122770, 2022.

\bibitem[Tan et~al.(2021)Tan, Xu, and Gao]{tan2021general}
Pengfei Tan, Tao Xu, and Feng Gao.
\newblock General coordinated active thermal control for parallel-connected inverters with switching frequency control.
\newblock In \emph{2021 IEEE 1st International Power Electronics and Application Symposium (PEAS)}, pages 1--6. IEEE, 2021.

\bibitem[Xu et~al.(2021)Xu, Zhao, and Liu]{xu2021near}
Deyu Xu, Junming Zhao, and Linhua Liu.
\newblock Near-field radiation assisted smart skin for spacecraft thermal control.
\newblock \emph{International Journal of Thermal Sciences}, 165:\penalty0 106934, 2021.

\bibitem[Zhang et~al.(2019)Zhang, Yang, Ma, Lin, and Yang]{zhang2019switched}
Xing Zhang, Qingqing Yang, Mingyao Ma, Zhengyu Lin, and Shuying Yang.
\newblock A switched reluctance motor torque ripple reduction strategy with deadbeat current control and active thermal management.
\newblock \emph{IEEE Transactions on Vehicular Technology}, 69\penalty0 (1):\penalty0 317--327, 2019.

\bibitem[Kuprat et~al.(2021)Kuprat, van~der Broeck, Andresen, Kalker, Liserre, and De~Doncker]{kuprat2021research}
Johannes Kuprat, Christoph~H van~der Broeck, Markus Andresen, Sven Kalker, Marco Liserre, and Rik~W De~Doncker.
\newblock Research on active thermal control: Actual status and future trends.
\newblock \emph{IEEE Journal of Emerging and Selected Topics in Power Electronics}, 9\penalty0 (6):\penalty0 6494--6506, 2021.

\bibitem[Miao et~al.(2021)Miao, Zhong, Zhao, and Zhao]{miao2021spacecraft}
Jianyin Miao, Qi~Zhong, Qiwei Zhao, and Xin Zhao.
\newblock \emph{Spacecraft thermal control technologies}.
\newblock Springer, 2021.

\bibitem[Taira et~al.(2017)Taira, Brunton, Dawson, Rowley, Colonius, McKeon, Schmidt, Gordeyev, Theofilis, and Ukeiley]{taira2017modal}
Kunihiko Taira, Steven~L Brunton, Scott~TM Dawson, Clarence~W Rowley, Tim Colonius, Beverley~J McKeon, Oliver~T Schmidt, Stanislav Gordeyev, Vassilios Theofilis, and Lawrence~S Ukeiley.
\newblock Modal analysis of fluid flows: An overview.
\newblock \emph{Aiaa Journal}, 55\penalty0 (12):\penalty0 4013--4041, 2017.

\bibitem[Brenner et~al.(2019)Brenner, Eldredge, and Freund]{brenner2019perspective}
MP~Brenner, JD~Eldredge, and JB~Freund.
\newblock Perspective on machine learning for advancing fluid mechanics.
\newblock \emph{Physical Review Fluids}, 4\penalty0 (10):\penalty0 100501, 2019.

\bibitem[Brunton et~al.(2020)Brunton, Noack, and Koumoutsakos]{brunton2020machine}
Steven~L Brunton, Bernd~R Noack, and Petros Koumoutsakos.
\newblock Machine learning for fluid mechanics.
\newblock \emph{Annual review of fluid mechanics}, 52:\penalty0 477--508, 2020.

\bibitem[Taira et~al.(2020)Taira, Hemati, Brunton, Sun, Duraisamy, Bagheri, Dawson, and Yeh]{taira2020modal}
Kunihiko Taira, Maziar~S Hemati, Steven~L Brunton, Yiyang Sun, Karthik Duraisamy, Shervin Bagheri, Scott~TM Dawson, and Chi-An Yeh.
\newblock Modal analysis of fluid flows: Applications and outlook.
\newblock \emph{AIAA journal}, 58\penalty0 (3):\penalty0 998--1022, 2020.

\bibitem[Bai et~al.(2017)Bai, Brunton, Brunton, Kutz, Kaiser, Spohn, and Noack]{bai2017data}
Zhe Bai, Steven~L Brunton, Bingni~W Brunton, J~Nathan Kutz, Eurika Kaiser, Andreas Spohn, and Bernd~R Noack.
\newblock \emph{Data-driven methods in fluid dynamics: Sparse classification from experimental data}.
\newblock Springer, 2017.

\bibitem[San et~al.(2019)San, Maulik, and Ahmed]{san2019artificial}
Omer San, Romit Maulik, and Mansoor Ahmed.
\newblock An artificial neural network framework for reduced order modeling of transient flows.
\newblock \emph{Communications in Nonlinear Science and Numerical Simulation}, 77:\penalty0 271--287, 2019.

\bibitem[Peng et~al.(2020)Peng, Chen, Aubry, Chen, and Wu]{peng2020time}
Jiang-Zhou Peng, Siheng Chen, Nadine Aubry, Zhi-Hua Chen, and Wei-Tao Wu.
\newblock Time-variant prediction of flow over an airfoil using deep neural network.
\newblock \emph{Physics of Fluids}, 32\penalty0 (12), 2020.

\bibitem[Peng et~al.(2021)Peng, Aubry, Zhu, Chen, and Wu]{peng2021geometry}
Jiang-Zhou Peng, Nadine Aubry, Shiquan Zhu, Zhihua Chen, and Wei-Tao Wu.
\newblock Geometry and boundary condition adaptive data-driven model of fluid flow based on deep convolutional neural networks.
\newblock \emph{Physics of Fluids}, 33\penalty0 (12), 2021.

\bibitem[Li et~al.(2020)Li, Zhang, Martins, and Shu]{li2020efficient}
Jichao Li, Mengqi Zhang, Joaquim~RRA Martins, and Chang Shu.
\newblock Efficient aerodynamic shape optimization with deep-learning-based geometric filtering.
\newblock \emph{AIAA journal}, 58\penalty0 (10):\penalty0 4243--4259, 2020.

\bibitem[Pawar et~al.(2020)Pawar, San, Rasheed, and Vedula]{pawar2020priori}
Suraj Pawar, Omer San, Adil Rasheed, and Prakash Vedula.
\newblock A priori analysis on deep learning of subgrid-scale parameterizations for kraichnan turbulence.
\newblock \emph{Theoretical and Computational Fluid Dynamics}, 34:\penalty0 429--455, 2020.

\bibitem[Vu et~al.(2021)Vu, Gulati, Vogel, Grunwald, and Bergs]{vu2021machine}
Anh~Tuan Vu, Shrey Gulati, Paul-Alexander Vogel, Tim Grunwald, and Thomas Bergs.
\newblock Machine learning-based predictive modeling of contact heat transfer.
\newblock \emph{International Journal of Heat and Mass Transfer}, 174:\penalty0 121300, 2021.

\bibitem[Swartz et~al.(2021)Swartz, Wu, Zhou, and Hao]{swartz2021machine}
Brandon Swartz, Lang Wu, Qiang Zhou, and Qing Hao.
\newblock Machine learning predictions of critical heat fluxes for pillar-modified surfaces.
\newblock \emph{International Journal of Heat and Mass Transfer}, 180:\penalty0 121744, 2021.

\bibitem[Wu et~al.(2019)Wu, Fang, and Xu]{wu2019predicting}
Yen-Ju Wu, Lei Fang, and Yibin Xu.
\newblock Predicting interfacial thermal resistance by machine learning.
\newblock \emph{npj Computational Materials}, 5\penalty0 (1):\penalty0 56, 2019.

\bibitem[Zou et~al.(2020)Zou, Yu, and Ergan]{zou2020towards}
Zhengbo Zou, Xinran Yu, and Semiha Ergan.
\newblock Towards optimal control of air handling units using deep reinforcement learning and recurrent neural network.
\newblock \emph{Building and Environment}, 168:\penalty0 106535, 2020.

\bibitem[Gao et~al.(2020)Gao, Li, and Wen]{gao2020deepcomfort}
Guanyu Gao, Jie Li, and Yonggang Wen.
\newblock Deepcomfort: Energy-efficient thermal comfort control in buildings via reinforcement learning.
\newblock \emph{IEEE Internet of Things Journal}, 7\penalty0 (9):\penalty0 8472--8484, 2020.

\bibitem[Wang and Hong(2020)]{wang2020reinforcement}
Zhe Wang and Tianzhen Hong.
\newblock Reinforcement learning for building controls: The opportunities and challenges.
\newblock \emph{Applied Energy}, 269:\penalty0 115036, 2020.

\bibitem[Xiong et~al.(2022)Xiong, Guo, and Tian]{xiong2022application}
Yan Xiong, Liang Guo, and Defu Tian.
\newblock Application of deep reinforcement learning to thermal control of space telescope.
\newblock \emph{Journal of Thermal Science and Engineering Applications}, 14\penalty0 (1):\penalty0 011011, 2022.

\bibitem[Rabault et~al.(2019)Rabault, Kuchta, Jensen, R{\'e}glade, and Cerardi]{rabault2019artificial}
Jean Rabault, Miroslav Kuchta, Atle Jensen, Ulysse R{\'e}glade, and Nicolas Cerardi.
\newblock Artificial neural networks trained through deep reinforcement learning discover control strategies for active flow control.
\newblock \emph{Journal of fluid mechanics}, 865:\penalty0 281--302, 2019.

\bibitem[Tang et~al.(2020)Tang, Rabault, Kuhnle, Wang, and Wang]{tang2020robust}
Hongwei Tang, Jean Rabault, Alexander Kuhnle, Yan Wang, and Tongguang Wang.
\newblock Robust active flow control over a range of reynolds numbers using an artificial neural network trained through deep reinforcement learning.
\newblock \emph{Physics of Fluids}, 32\penalty0 (5), 2020.

\bibitem[Xu et~al.(2020)Xu, Zhang, Deng, and Rabault]{xu2020active}
Hui Xu, Wei Zhang, Jian Deng, and Jean Rabault.
\newblock Active flow control with rotating cylinders by an artificial neural network trained by deep reinforcement learning.
\newblock \emph{Journal of Hydrodynamics}, 32\penalty0 (2):\penalty0 254--258, 2020.

\bibitem[Li and Zhang(2022)]{li2022reinforcement}
Jichao Li and Mengqi Zhang.
\newblock Reinforcement-learning-based control of confined cylinder wakes with stability analyses.
\newblock \emph{Journal of Fluid Mechanics}, 932:\penalty0 A44, 2022.

\bibitem[Ren et~al.(2021)Ren, Rabault, and Tang]{ren2021applying}
Feng Ren, Jean Rabault, and Hui Tang.
\newblock Applying deep reinforcement learning to active flow control in weakly turbulent conditions.
\newblock \emph{Physics of Fluids}, 33\penalty0 (3), 2021.

\bibitem[Shimomura et~al.(2020)Shimomura, Sekimoto, Oyama, Fujii, and Nishida]{shimomura2020experimental}
Satoshi Shimomura, Satoshi Sekimoto, Akira Oyama, Kozo Fujii, and Hiroyuki Nishida.
\newblock Experimental study on application of distributed deep reinforcement learning to closed-loop flow separation control over an airfoil.
\newblock In \emph{AIAA Scitech 2020 Forum}, page 0579, 2020.

\bibitem[Fan et~al.(2020)Fan, Yang, Wang, Triantafyllou, and Karniadakis]{fan2020reinforcement}
Dixia Fan, Liu Yang, Zhicheng Wang, Michael~S Triantafyllou, and George~Em Karniadakis.
\newblock Reinforcement learning for bluff body active flow control in experiments and simulations.
\newblock \emph{Proceedings of the National Academy of Sciences}, 117\penalty0 (42):\penalty0 26091--26098, 2020.

\bibitem[Wang et~al.(2022)Wang, Mei, Aubry, Chen, Wu, and Wu]{wang2022deep}
Yi-Zhe Wang, Yu-Fei Mei, Nadine Aubry, Zhihua Chen, Peng Wu, and Wei-Tao Wu.
\newblock Deep reinforcement learning based synthetic jet control on disturbed flow over airfoil.
\newblock \emph{Physics of Fluids}, 34\penalty0 (3), 2022.

\bibitem[Mei et~al.(2021)Mei, Zheng, Aubry, Li, Wu, and Liu]{mei2021active}
Yu-Fei Mei, Chun Zheng, Nadine Aubry, Meng-Ge Li, Wei-Tao Wu, and Xianglei Liu.
\newblock Active control for enhancing vortex induced vibration of a circular cylinder based on deep reinforcement learning.
\newblock \emph{Physics of Fluids}, 33\penalty0 (10), 2021.

\bibitem[Doll et~al.(2012)Doll, Simon, and Daw]{doll2012ubiquity}
Bradley~B Doll, Dylan~A Simon, and Nathaniel~D Daw.
\newblock The ubiquity of model-based reinforcement learning.
\newblock \emph{Current opinion in neurobiology}, 22\penalty0 (6):\penalty0 1075--1081, 2012.

\bibitem[Kaiser et~al.(2019)Kaiser, Babaeizadeh, Milos, Osinski, Campbell, Czechowski, Erhan, Finn, Kozakowski, Levine, et~al.]{kaiser2019model}
Lukasz Kaiser, Mohammad Babaeizadeh, Piotr Milos, Blazej Osinski, Roy~H Campbell, Konrad Czechowski, Dumitru Erhan, Chelsea Finn, Piotr Kozakowski, Sergey Levine, et~al.
\newblock Model-based reinforcement learning for atari.
\newblock \emph{arXiv preprint arXiv:1903.00374}, 2019.

\bibitem[Polydoros and Nalpantidis(2017)]{polydoros2017survey}
Athanasios~S Polydoros and Lazaros Nalpantidis.
\newblock Survey of model-based reinforcement learning: Applications on robotics.
\newblock \emph{Journal of Intelligent \& Robotic Systems}, 86\penalty0 (2):\penalty0 153--173, 2017.

\bibitem[{\c{C}}al{\i}{\c{s}}{\i}r and Pehlivano{\u{g}}lu(2019)]{ccalicsir2019model}
Sinan {\c{C}}al{\i}{\c{s}}{\i}r and Meltem~Kurt Pehlivano{\u{g}}lu.
\newblock Model-free reinforcement learning algorithms: A survey.
\newblock In \emph{2019 27th signal processing and communications applications conference (SIU)}, pages 1--4. IEEE, 2019.

\bibitem[ANSYS(2019)]{ansys2019ansys}
Inc. ANSYS.
\newblock Ansys fluent user’s guide, version 2019 r3, 2019.

\bibitem[Jang et~al.(2019)Jang, Kim, Harerimana, and Kim]{jang2019q}
Beakcheol Jang, Myeonghwi Kim, Gaspard Harerimana, and Jong~Wook Kim.
\newblock Q-learning algorithms: A comprehensive classification and applications.
\newblock \emph{IEEE access}, 7:\penalty0 133653--133667, 2019.

\bibitem[Mnih et~al.(2013)Mnih, Kavukcuoglu, Silver, Graves, Antonoglou, Wierstra, and Riedmiller]{mnih2013playing}
Volodymyr Mnih, Koray Kavukcuoglu, David Silver, Alex Graves, Ioannis Antonoglou, Daan Wierstra, and Martin Riedmiller.
\newblock Playing atari with deep reinforcement learning.
\newblock \emph{arXiv preprint arXiv:1312.5602}, 2013.

\bibitem[Mnih et~al.(2015)Mnih, Kavukcuoglu, Silver, Rusu, Veness, Bellemare, Graves, Riedmiller, Fidjeland, Ostrovski, et~al.]{mnih2015human}
Volodymyr Mnih, Koray Kavukcuoglu, David Silver, Andrei~A Rusu, Joel Veness, Marc~G Bellemare, Alex Graves, Martin Riedmiller, Andreas~K Fidjeland, Georg Ostrovski, et~al.
\newblock Human-level control through deep reinforcement learning.
\newblock \emph{nature}, 518\penalty0 (7540):\penalty0 529--533, 2015.

\bibitem[Kingma and Ba(2014)]{kingma2014adam}
Diederik~P Kingma and Jimmy Ba.
\newblock Adam: A method for stochastic optimization.
\newblock \emph{arXiv preprint arXiv:1412.6980}, 2014.

\bibitem[Hessel et~al.(2018)Hessel, Modayil, Van~Hasselt, Schaul, Ostrovski, Dabney, Horgan, Piot, Azar, and Silver]{hessel2018rainbow}
Matteo Hessel, Joseph Modayil, Hado Van~Hasselt, Tom Schaul, Georg Ostrovski, Will Dabney, Dan Horgan, Bilal Piot, Mohammad Azar, and David Silver.
\newblock Rainbow: Combining improvements in deep reinforcement learning.
\newblock In \emph{Proceedings of the AAAI conference on artificial intelligence}, volume~32, 2018.

\bibitem[Arulkumaran et~al.(2017)Arulkumaran, Deisenroth, Brundage, and Bharath]{arulkumaran2017brief}
Kai Arulkumaran, Marc~Peter Deisenroth, Miles Brundage, and Anil~Anthony Bharath.
\newblock A brief survey of deep reinforcement learning.
\newblock \emph{arXiv preprint arXiv:1708.05866}, 2017.

\bibitem[Sewak(2019)]{sewak2019deep}
Mohit Sewak.
\newblock \emph{Deep reinforcement learning}.
\newblock Springer, 2019.

\bibitem[Wan et~al.(2020)Wan, Gao, Hu, and Wu]{wan2020robust}
Kaifang Wan, Xiaoguang Gao, Zijian Hu, and Gaofeng Wu.
\newblock Robust motion control for uav in dynamic uncertain environments using deep reinforcement learning.
\newblock \emph{Remote sensing}, 12\penalty0 (4):\penalty0 640, 2020.

\bibitem[Haarnoja et~al.(2018)Haarnoja, Pong, Zhou, Dalal, Abbeel, and Levine]{haarnoja2018composable}
Tuomas Haarnoja, Vitchyr Pong, Aurick Zhou, Murtaza Dalal, Pieter Abbeel, and Sergey Levine.
\newblock Composable deep reinforcement learning for robotic manipulation.
\newblock In \emph{2018 IEEE international conference on robotics and automation (ICRA)}, pages 6244--6251. IEEE, 2018.

\bibitem[Hasselt(2010)]{hasselt2010double}
Hado Hasselt.
\newblock Double q-learning.
\newblock \emph{Advances in neural information processing systems}, 23, 2010.

\bibitem[Van~Hasselt et~al.(2016)Van~Hasselt, Guez, and Silver]{van2016deep}
Hado Van~Hasselt, Arthur Guez, and David Silver.
\newblock Deep reinforcement learning with double q-learning.
\newblock In \emph{Proceedings of the AAAI conference on artificial intelligence}, volume~30, 2016.

\bibitem[Wang et~al.(2016)Wang, Schaul, Hessel, Hasselt, Lanctot, and Freitas]{wang2016dueling}
Ziyu Wang, Tom Schaul, Matteo Hessel, Hado Hasselt, Marc Lanctot, and Nando Freitas.
\newblock Dueling network architectures for deep reinforcement learning.
\newblock In \emph{International conference on machine learning}, pages 1995--2003. PMLR, 2016.

\bibitem[Tavakoli et~al.(2018)Tavakoli, Pardo, and Kormushev]{tavakoli2018action}
Arash Tavakoli, Fabio Pardo, and Petar Kormushev.
\newblock Action branching architectures for deep reinforcement learning.
\newblock In \emph{Proceedings of the aaai conference on artificial intelligence}, volume~32, 2018.

\bibitem[Sodjavi et~al.(2015)Sodjavi, Montagn{\'e}, Bragan{\c{c}}a, Meslem, Bode, and Kristiawan]{sodjavi2015impinging}
Kodjovi Sodjavi, Brice Montagn{\'e}, Pierre Bragan{\c{c}}a, Amina Meslem, Florin Bode, and Magdalena Kristiawan.
\newblock Impinging cross-shaped submerged jet on a flat plate: a comparison of plane and hemispherical orifice nozzles.
\newblock \emph{Meccanica}, 50:\penalty0 2927--2947, 2015.

\bibitem[Bunker et~al.(2014)Bunker, Dees, and Palafox]{bunker2014impingement}
Ronald~S Bunker, Jason~E Dees, and Pepe Palafox.
\newblock Impingement cooling in gas turbines: Design, applications, and limitations.
\newblock \emph{Impingement Jet Cooling in Gas Turbines}, 25\penalty0 (1), 2014.

\bibitem[Wang et~al.(2020)Wang, Guo, Xiong, and Wang]{wang2020review}
Ji-Xiang Wang, Wei Guo, Kai Xiong, and Sheng-Nan Wang.
\newblock Review of aerospace-oriented spray cooling technology.
\newblock \emph{Progress in Aerospace Sciences}, 116:\penalty0 100635, 2020.

\bibitem[Jena et~al.(2022)Jena, Mishra, and Sahoo]{jena2022comparative}
Mrutyunjay Jena, Purna~Chandra Mishra, and Sudhansu~Sekhar Sahoo.
\newblock Comparative performance of jet and spray impingement cooling in steel strip run-out table: experimental results.
\newblock \emph{Australian Journal of Mechanical Engineering}, pages 1--14, 2022.

\bibitem[Wei et~al.(2018)Wei, Wicke, Freelan, and Luke]{wei2018multiagent}
Ermo Wei, Drew Wicke, David Freelan, and Sean Luke.
\newblock Multiagent soft q-learning.
\newblock \emph{arXiv preprint arXiv:1804.09817}, 2018.

\end{thebibliography}


\end{document}